# Hypocoercivity for non-linear infinite-dimensional degenerate stochastic differential equations

Benedikt Eisenhuth[1] · Martin Grothaus[1]

## Abstract

The aim of this article is to construct solutions to second order in time stochastic partial differential equations and to show hypocoercivity of the corresponding transition semigroups. More generally, we analyze non-linear infinite-dimensional degenerate stochastic differential equations in terms of their infinitesimal generators. In the first part of this article we use resolvent methods developed by Beznea, Boboc and Röckner to construct diffusion processes with infinite lifetime and explicit invariant measures. The processes provide weak solutions to infinite-dimensional Langevin dynamics. The second part deals with a general abstract Hilbert space hypocoercivity method, developed by Dolbeaut, Mouhot and Schmeiser. In order to treat stochastic (partial) differential equations, Grothaus and Stilgenbauer translated these concepts to the Kolmogorov backwards setting taking domain issues into account. To apply these concepts in the context of infinite-dimensional Langevin dynamics we use an essential m-dissipativity result for infinite-dimensional Ornstein–Uhlenbeck operators, perturbed by the gradient of a potential. We allow unbounded diffusion operators as coefficients and apply corresponding regularity estimates. Furthermore, essential m-dissipativity of a non-sectorial Kolmogorov backward operator associated to the dynamic is needed. Poincaré inequalities for measures with densities w.r.t. infinite-dimensional non-degenerate Gaussian measures are studied. Deriving a stochastic representation of the semigroup generated by the Kolmogorov backward operator as the transition semigroup of a diffusion process enables us to show an $L^2$-exponential ergodicity result for the weak solution. Finally, we apply our results to explicit infinite-dimensional degenerate diffusion equations.

✉ Benedikt Eisenhuth
  eisenhuth@mathematik.uni-kl.de

  Martin Grothaus
  grothaus@mathematik.uni-kl.de

[1] Department of Mathematics, RPTU Kaiserslautern-Landau, PO Box 3049, 67653 Kaiserslautern, Germany

## 1 Introduction

Let $(U, (\cdot, \cdot)_U)$ and $(V, (\cdot, \cdot)_V)$ be two real separable Hilbert spaces. We equip $W = U \times V$ with the inner product $(\cdot, \cdot)_W$ defined by

$$((u_1, v_1), (u_2, v_2))_W = (u_1, u_2)_U + (v_1, v_2)_V, \quad (u_1, v_1), (u_2, v_2) \in W.$$

Denote by $\mathcal{B}(U)$ and $\mathcal{B}(V)$ the Borel $\sigma$-algebra on $U$ and $V$, on which we consider centered non-degenerate Gaussian measures $\mu_1$ and $\mu_2$, respectively. The measures are uniquely determined by their covariance operators $Q_1 \in \mathcal{L}(U)$ and $Q_2 \in \mathcal{L}(V)$, respectively. Furthermore, we fix bounded linear operators $K_{12} \in \mathcal{L}(U, V)$, $K_{21} \in \mathcal{L}(V, U)$, a symmetric positive bounded linear operator $K_{22} \in \mathcal{L}(V)$ and a measurable potential $\Phi: U \to (-\infty, \infty]$. We assume that $\Phi$ is bounded from below and $\int_U e^{-\Phi} d\mu_1 = 1$. For such $\Phi$ we consider the probability measure $\mu_1^\Phi = e^{-\Phi} \mu_1$ on $(U, \mathcal{B}(U))$ and the product measure

$$\mu^\Phi = \mu_1^\Phi \otimes \mu_2,$$

on $(W, \mathcal{B}(W))$. In [16] we analyzed the infinite-dimensional Langevin operator $(L_\Phi, \mathcal{F}C_b^\infty(B_W))$ in $L^2(W, \mu^\Phi, \mathbb{R})$ defined by

$$L_\Phi = S - A_\Phi,$$

where $S$ and $A_\Phi$ applied to $f \in \mathcal{F}C_b^\infty(B_W)$ are given via

$$Sf = \mathrm{tr}[K_{22} D_2^2 f] - (v, Q_2^{-1} K_{22} D_2 f)_V \quad \text{and}$$
$$A_\Phi f = (u, Q_1^{-1} K_{21} D_2 f)_U + (D\Phi(u), K_{21} D_2 f)_U - (v, Q_2^{-1} K_{12} D_1 f)_V.$$

Here $\mathcal{F}C_b^\infty(B_W)$ denotes the space of finitely based smooth (infinitely often differentiable with bounded derivatives) cylinder functions on $W$ and $D_1$, $D_2$ and $D_2^2$ the first and second order Fréchet derivatives in the first and second component, respectively, see Remark 2.1.

An application of the Itô-formula shows that $(L_\Phi, \mathcal{F}C_b^\infty(B_W))$ corresponds to the infinite-dimensional non-linear stochastic differential equation, heuristically given by

$$\begin{aligned} dU_t &= K_{21} Q_2^{-1} V_t dt \\ dV_t &= -(K_{22} Q_2^{-1} V_t + K_{12} Q_1^{-1} U_t + K_{12} D\Phi(U_t)) dt + \sqrt{2 K_{22}} dW_t. \end{aligned} \quad (1)$$

In this article $(W_t)_{t\geq 0}$ denotes a cylindrical Brownian motion with values in $V$. Equation (1) is an infinite-dimensional Langevin equation. We emphasize its degenerate structure, i.e. the noise only appears in the second component. This corresponds to the fact, that the antisymmetric part $(A_\Phi, \mathcal{F}C_b^\infty(B_W))$ of $(L_\Phi, \mathcal{F}C_b^\infty(B_W))$ in $L^2(W, \mu^\Phi, \mathbb{R})$ contains first order differential operators in the first component and the symmetric part $(S, \mathcal{F}C_b^\infty(B_W))$ only differential operators in the second component. In particular the operator $(L_\Phi, \mathcal{F}C_b^\infty(B_W))$ is non-sectorial.

In order to show exponential convergence to equilibrium of such non-coercive evolution equations, Cédric Villani developed the concept of hypocoercivity, see [27]. Abstract hypocoercivity concepts for a quantitative description of convergence rates were introduced in [13]. In order to treat stochastic (partial) differential equations, Grothaus and Stilgenbauer translated these concepts to the Kolmogorov backwards setting taking domain issues into account, see [17, 18]. In [20] these concepts were further generalized to the case where only a weak Poincaré inequality is needed. In this case one obtains (sub-)exponential convergence rates. Ergodicity and the rate of convergence to equilibrium of finite-dimensional Langevin dynamics with singular potentials were studied in [6, 19]. There the authors used generalized Dirichlet forms and martingale techniques. Singular potentials could be treated also by Lyapunov techniques, see [2, 5, 21]. Using coupling methods, Langevin dynamics and their quantitative contraction rates were studied in [14]. Recently, corresponding dynamics and their hypocoercive behavior were studied on abstract smooth manifolds, see [24] or with multiplicative noise, see [3]. In [28] infinite-dimensional stochastic Hamiltonian systems were treated. The results obtained there are complementary to ours. There the author used coupling methods and dimension-free Harnack inequalities in order to prove hypercontractivity in a degenerate situation.

The structure of this article is as follows. Basic notations are stated in Sect. 2. In Sect. 3 we state the preliminary results needed in this article and recently obtained in [16]. E.g. an essential m-dissipativity result for infinite-dimensional Ornstein–Uhlenbeck operators, perturbed by the gradient of a potential, with possible unbounded diffusion operators as coefficients, see Theorem 3.11. Furthermore, we state corresponding first and second order regularity estimates, see again Theorem 3.11. Considering bounded diffusion operators as coefficients, similar results and estimates have been used in [10, Section 4] to study stochastic quantization problems. Non-degenerate stochastic reaction–diffusion equations have been analyzed in [9, Section 5]. In Theorem 3.6, we quote the statement of essential m-dissipativity of $(L_\Phi, \mathcal{F}C_b^\infty(B_W))$ in $L^2(W, \mu^\Phi, \mathbb{R})$, derived in [16].

In Sect. 4 we construct weak solutions, via the martingale problem, to the infinite-dimensional Langevin equation. The main strategy to achieve this is discussed in the article [4] from Lucian Beznea, Nicu Boboc and Michael Röckner. They provide a general strategy to construct martingale solutions in terms of a $\mu^\Phi$-standard right process

$$X = (\Omega, \mathcal{F}, (\mathcal{F}_t)_{t\geq 0}, (X_t)_{t\geq 0}, (\theta_t)_{t\geq 0}, (\mathbb{P}^w)_{w\in W}),$$

which is càdlàg. The transition semigroup of the process $X$ extended to $L^2(W, \mu^\Phi, \mathbb{R})$ coincides with the strongly continuous contraction semigroup associated to $(L_\Phi, \mathcal{F}C_b^\infty(B_W))$. The application of this abstract existence result requires that $K_{22}$ is of trace class, compare Assumption 3.8. Similar to [7, Theorem 6.4.2.], we show that $(X_t)_{t\geq 0}$ has infinite lifetime and weakly continuous paths. In Theorem 4.7 we discuss the existence of weak solutions. We construct a cylindrical Brownian motion $(W_t)_{t\geq 0}$ with values in $V$ s.t. $\mathbb{P}_{\mu^\Phi}$-a.s. it holds for all $\vartheta \in D(Q_2^{-1}K_{12})$, $\theta \in D(Q_1^{-1}K_{21}) \cap D(Q_2^{-1}K_{22})$ and $t \in [0, \infty)$

$$(U_t, \vartheta)_U = (U_0, \vartheta)_U + \int_{[0,t]} (V_s, Q_2^{-1}K_{12}\vartheta)_V \, d\lambda(s),$$

$$(V_t, \theta)_V = (V_0, \theta)_V - \int_{[0,t]} (V_s, Q_2^{-1}K_{22}\theta)_V + (U_s, Q_1^{-1}K_{21}\theta)_U$$
$$+ (D\Phi(U_s), K_{21}\theta)_U \, d\lambda(s) + (\sqrt{2K_{22}}W_t, \theta)_V.$$

Above, $(X_t)_{t\geq 0} = ((U_t, V_t))_{t\geq 0}$ and $\mathbb{P}_{\mu^\Phi} = \int_W \mathbb{P}^w d\mu^\Phi(w)$.

The general abstract Hilbert space hypocoercivity method from [17] is described in Sect. 5. In Theorem 5.2 we show, on an abstract level, how to explicitly compute the exponential convergence rate. Section 6 contains the application of this method in the infinite-dimensional Langevin setting. There, the essential m-dissipativity and the regularity results proved in [16] are necessary. Moreover, the Poincaré type inequalities from Assumption 3.10 are essential. Indeed, they are used in Theorem 6.10, to explicitly compute, for all $\theta_1 \in (1, \infty)$, a constant $\theta_2 \in (0, \infty)$ such that for each $g \in L^2(W, \mu^\Phi, \mathbb{R})$ and all $t \in [0, \infty)$ we have

$$\|T_t g - (g, 1)_{L^2(\mu^\Phi)}\|_{L^2(\mu^\Phi)} \leq \theta_1 e^{-\theta_2 t} \|g - (g, 1)_{L^2(\mu^\Phi)}\|_{L^2(\mu^\Phi)}.$$

Above, $(T_t)_{t\geq 0}$ is the strongly continuous contraction semigroup generated by the closure of $(L_\Phi, \mathcal{F}C_b^\infty(B_W))$. Up to our knowledge this is the first time, this abstract hypocoercivity method is applied to show hypocoercivity of an infinite-dimensional degenerate stochastic differential equation. We would like to stress that the hypocoercivity result for the strongly continuous contraction semigroup $(T_t)_{t\geq 0}$ also applies, if $K_{22}$ is not a trace class operator. The latter is needed to construct the stochastic process $X$. At the end of Sect. 6 we describe, assuming $K_{22}$ is a trace class operator, the connection between the transition semigroup of $X$ and $(T_t)_{t\geq 0}$. Moreover, we provide an $L^2$-exponential ergodicity result for the weak solution of the infinite-dimensional Langevin equation. In the last section, we present explicit examples to which the results from the previous sections are applicable. We would like to emphasize that our results can treat equations, beyond the framework of degenerate semi-linear stochastic partial differential equations, as considered in [28]. In particular, we are able to analyze degenerate stochastic reaction–diffusion equations. The main results obtained in this article are summarized in the following list:

- We construct a $\mu^\Phi$-standard right process $X$ with infinite lifetime and weakly continuous paths, providing a weak solution, via the martingale problem, to the infinite-dimension Langevin equation, see Theorems 4.3 and 4.7.

- Hypocoercivity, i.e. exponential convergence to equilibrium of the strongly continuous semigroup generated by the infinite-dimensional Langevin operator $(L_\Phi, \mathcal{F}C_b^\infty(B_W))$ in the Hilbert space $L^2(W, \mu^\Phi, \mathbb{R})$ with explicitly computable rate of convergence is proved in Theorem 6.10.
- We derive a stochastic representation of the semigroup $(T_t)_{t\geq 0}$ in terms of the transition semigroup of the $\mu^\Phi$-standard right process $X$. Corollary 6.12 contains an $L^2$-exponential ergodicity result, corresponding to the hypocoercivity of the semigroup generated by the infinite-dimensional Langevin operator $(L_\Phi, \mathcal{F}C_b^\infty(B_W))$.

## 2 Notations

Let $U$ and $V$ be two real separable Hilbert spaces with inner products $(\cdot, \cdot)_U$ and $(\cdot, \cdot)_V$, respectively. In particular there exist a orthonormal basis $B_U = (d_i)_{i \in \mathbb{N}}$ and $B_V = (e_i)_{i \in \mathbb{N}}$ of $U$ and $V$, respectively. Let $W = U \times V$ be the real separable Hilbert space, equipped with the canonical inner product considered in the introduction. $B_W = \{(d_n, 0) \mid n \in \mathbb{N}\} \cup \{(0, e_n) \mid n \in \mathbb{N}\}$ is an orthonormal basis of $W$.

The set of all linear bounded operators from $U$ to $U$ and from $U$ to $V$ are denoted by $\mathcal{L}(U)$ and $\mathcal{L}(U, V)$, respectively. The adjoint of an operator $K \in \mathcal{L}(U, V)$ is denoted by $K^* \in \mathcal{L}(V, U)$. By $\mathcal{L}^+(U)$ we shall denote the subset of $\mathcal{L}(U)$ consisting of all positive symmetric operators. The subset in $\mathcal{L}^+(U)$ of trace class operators is denoted by $\mathcal{L}_1^+(U)$. The Borel $\sigma$-algebra on a Hilbert space $(Y, (\cdot, \cdot)_Y)$ is denoted by $\mathcal{B}(Y)$. For a given measure space $(\Omega, \mathcal{A}, m)$ we denote by $L^p(\Omega, m, Y)$, $p \in (0, \infty]$, the Banach space of equivalence classes of $\mathcal{A}\text{-}\mathcal{B}(Y)$ measurable and $p$-integrable functions. The corresponding norm is denoted by $\|\cdot\|_{L^p(\Omega, m, Y)}$. If $p = 2$ the norm is induced by an inner product denoted by $(\cdot, \cdot)_{L^2(\Omega, m, Y)}$. In case $(\Omega, \mathcal{A})$ and and $Y$ is clear from the context we also write $L^2(m)$ instead of $L^2(\Omega, m, \mathbb{R}^n)$. By $\lambda$ we denote the Lebesgue measure on $(\mathbb{R}, \mathcal{B}(\mathbb{R}))$.

**Remark 2.1** Given a Fréchet differentiable function $f : U \to \mathbb{R}$. For $u \in U$ we denote by $Df(u) \in U$ the Fréchet derivative of $f$ in $u$. In this article we identify the Fréchet derivative with the gradient. Analogously, for a two times Fréchet differentiable function $f : U \to \mathbb{R}$, we identify the second order Fréchet derivative $D^2 f(u)$ in $u \in U$ with a bounded linear map in $L(U)$. For $i, j \in \mathbb{N}$ we denote by $\partial_i f(u) = (Df(u), d_i)_U$ the partial derivative in the direction of $d_i$ and by $\partial_{ij} f(u) = (D^2 f(u) d_i, d_j)_U$ the second order partial derivative in the direction of $d_i$ and $d_j$. In particular we have

$$Df = \sum_{i=1}^\infty (Df, d_i)_U d_i = \sum_{i=1}^\infty \partial_i f d_i.$$

Concerning derivatives of sufficient smooth functions $f : W \to \mathbb{R}$ we set $D_1 f = \sum_{i=1}^\infty (Df, (d_i, 0))_W d_i \in U$ and $D_2 f = \sum_{i=1}^\infty (Df, (0, e_i))_W e_i \in V$ as well as $\partial_{i,1} f = (D_1 f, d_i)_U$ and $\partial_{i,2} f = (D_2 f, e_i)_V$. In particular we have

$$Df = \sum_{i=1}^\infty (Df, (d_i, 0))_W (d_i, 0) + \sum_{i=1}^\infty (Df, (0, e_i))_W (0, e_i) = (D_1 f, D_2 f).$$

Analogously we define $D_1^2 f$, $D_2^2 f$ as well as $\partial_{ij,1} f$ and $\partial_{ij,2} f$.

The orthogonal projection from $U$ to $B_U^n = \text{span}\{d_1, ..., d_n\}$ is denoted by $\overline{P}_n$ and the corresponding coordinate map by $P_n$, i.e. we have for all $u \in U$,

$$\overline{P}_n(u) = \sum_{i=1}^n (u, d_i)_U d_i \quad \text{and} \quad P_n(u) = \big((u, d_1)_U, ..., (u, d_n)_U\big).$$

To avoid an overload of notation we use the same notation for the projection to $B_V^n = \text{span}\{e_1, ..., e_n\}$. The spaces of finitely based smooth (infinitely often differentiable with bounded derivatives) cylinder functions w.r.t. $B_U$ and $B_W$ are denoted by $\mathcal{F}C_b^\infty(B_U)$ and $\mathcal{F}C_b^\infty(B_W)$, respectively. As in [16] these spaces are defined by

$$\mathcal{F}C_b^\infty(B_U) = \bigcup_{n \in \mathbb{N}} \mathcal{F}C_b^\infty(B_U, n) \quad \text{and} \quad \mathcal{F}C_b^\infty(B_W) = \bigcup_{n \in \mathbb{N}} \mathcal{F}C_b^\infty(B_W, n),$$

where $\mathcal{F}C_b^\infty(B_U, n)$ and $\mathcal{F}C_b^\infty(B_W, n)$ denote the spaces of finitely based smooth cylinder functions only dependent on the first $n$ directions defined by

$$\mathcal{F}C_b^\infty(B_U, n) = \{U \ni u \mapsto \varphi(P_n(u)) \in \mathbb{R} \mid \varphi \in C_b^\infty(\mathbb{R}^n)\} \quad \text{and}$$
$$\mathcal{F}C_b^\infty(B_W, n) = \{W \ni (u, v) \mapsto \varphi(P_n(u), P_n(v)) \in \mathbb{R} \mid \varphi \in C_b^\infty(\mathbb{R}^n \times \mathbb{R}^n)\}.$$

Above, $C_b^\infty(\mathbb{R}^m)$ and $C_b^\infty(\mathbb{R}^m \times \mathbb{R}^m)$ denote the space of infinitively differentiable bounded functions from $\mathbb{R}^m$ and $\mathbb{R}^m \times \mathbb{R}^m$ to $\mathbb{R}$, respectively. By $C_0^\infty(\mathbb{R}^n)$ we denote the space of compactly supported smooth functions from $\mathbb{R}^n$ to $\mathbb{R}$.

## 3 Assumptions and preliminaries

Let $U$, $V$ and $W$ as above. We fix two centered non-degenerate Gaussian measure $\mu_1$ and $\mu_2$ with covariance operators $Q_1 \in \mathcal{L}_1^+(U)$ and $Q_2 \in \mathcal{L}_1^+(V)$ on $(U, \mathcal{B}(U))$ and $(V, \mathcal{B}(V))$, respectively. Let $B_U = (d_i)_{i \in \mathbb{N}}$ and $B_V = (e_i)_{i \in \mathbb{N}}$ be the orthonormal basis of eigenvectors with corresponding eigenvalues $(\lambda_i)_{i \in \mathbb{N}}$ and $(\nu_i)_{i \in \mathbb{N}}$ of $Q_1$ and $Q_2$, respectively. W.l.o.g. we assume that $(\lambda_i)_{i \in \mathbb{N}}$ and $(\nu_i)_{i \in \mathbb{N}}$ are decreasing to zero. On $(W, \mathcal{B}(W))$ we define the measure $\mu$ by

$$\mu = \mu_1 \otimes \mu_2.$$

In [16] it is shown, that the function spaces $\mathcal{F}C_b^\infty(B_U)$ and $\mathcal{F}C_b^\infty(B_W)$ are dense in $L^2(U, \mu_1, \mathbb{R})$ and $L^2(W, \mu, \mathbb{R})$, respectively.

**Remark 3.1** The domain of the closure of

$$D : \mathcal{F}C_b^\infty(B_U) \to L^2(U, \mu_1, U),$$

in $L^2(U, \mu_1, U)$ is denoted by $W^{1,2}(U, \mu_1, \mathbb{R})$, compare [16, Theorem 2.9] and [11, Chapter 9.2]. I.e. $W^{1,2}(U, \mu_1, \mathbb{R})$ is the first order Sobolev space w.r.t. the infinite-dimensional Gaussian measure $\mu_1$. The closure of $D$ is again denoted by $D$ and is referred to as the weak gradient in the following. By [8, Proposition 10.6] every bounded Fréchet differentiable function $f : U \to \mathbb{R}$ with bounded gradient is in $W^{1,2}(U, \mu_1, \mathbb{R})$ and the the classical gradient of $f$ coincides $\mu_1$-a.e. with the weak one.

We fix a measurable potential $\Phi : U \to (-\infty, \infty]$, which is bounded from below. Moreover, during the whole article we assume $\int_U e^{-\Phi} d\mu_1 = 1$ and $\Phi \in W^{1,2}(U, \mu_1, \mathbb{R})$. For such potentials $\Phi$, we consider the measure $\mu_1^\Phi = e^{-\Phi} \mu_1$ on $(U, \mathcal{B}(U))$. On $(W, \mathcal{B}(W))$ we define the measure $\mu^\Phi$ by

$$\mu^\Phi = \mu_1^\Phi \otimes \mu_2.$$

Note that $\mathcal{F}C_b^\infty(B_U)$ and $\mathcal{F}C_b^\infty(B_W)$ are also dense in $L^2(U, \mu_1^\Phi, \mathbb{R})$ and $L^2(W, \mu^\Phi, \mathbb{R})$, respectively.

Another well known approximation result providing a particular sequence of smooth functions in $\mathbb{R}^n$ is stated in the next lemma. Such sequences are needed several times during the article.

**Lemma 3.2** *Let $m \in \mathbb{N}$. There is some $\varphi \in C_0^\infty(\mathbb{R}^n)$ such that $0 \leq \varphi \leq 1$, $\varphi = 1$ on $B_1(0) = \{x \in \mathbb{R}^n \mid |x| < 1\}$ and $\varphi = 0$ outside $B_2(0)$ and a constant $c \in (0, \infty)$, independent of m, such that*

$$|\partial_i \varphi_m(z)| \leq \frac{c}{m}, \quad |\partial_{ij} \varphi_m(z)| \leq \frac{c}{m^2} \quad \text{for all} \quad z \in \mathbb{R}^n, \quad 1 \leq i, j \leq n,$$

*where we define*

$$\varphi_m(z) = \varphi(\frac{z}{m}) \quad \text{for each} \quad z \in \mathbb{R}^n.$$

*In particular, $0 \leq \varphi_m \leq 1$ and $\varphi_m = 1$ on $B_m(0)$ for all $m \in \mathbb{N}$. Moreover, $\varphi_m \to 1$ pointwisely on $\mathbb{R}^n$ and $D\varphi_m, D^2\varphi \to 0$ as $m \to \infty$, w.r.t. the sup norm $|\cdot|_\infty$.*

During the article we have to calculate certain Gaussian integrals as written down in the lemma below. The proof uses Isserlis formula from [22] and can be found in [16].

**Lemma 3.3** [16, Lemma 2.2] *For $l_1, l_2, l_3, l_4 \in U$ set $q_{ij} = (Q_1 l_i, l_j)_U$, $i, j \in \{1, 2, 3, 4\}$. Then it holds*

$$\int_U (u, l_1)_U (u, l_2)_U d\mu_1(u) = q_{12} \quad \text{and}$$

$$\int_U (u, l_1)_U (u, l_2)_U (u, l_3)_U (u, l_4)_U d\mu_1(u) = q_{12}q_{34} + q_{13}q_{24} + q_{14}q_{23}.$$

As already described in the introduction we define the infinite-dimensional Langevin operator $(L_\Phi, \mathcal{F}C_b^\infty(B_W))$ in $L^2(W, \mu^\Phi, \mathbb{R})$ by

$$L_\Phi = S - A_\Phi,$$

where for $f \in \mathcal{F}C_b^\infty(B_W)$, $Sf$ and $A_\Phi f$ are given via

$$Sf = \text{tr}[K_{22} D_2^2 f] - (v, Q_2^{-1} K_{22} D_2 f)_V \quad \text{and}$$
$$A_\Phi f = (u, Q_1^{-1} K_{21} D_2 f)_U + (D\Phi(u), K_{21} D_2 f)_U - (v, Q_2^{-1} K_{12} D_1 f)_V.$$

Throughout the whole article we assume the following:

**Assumption 3.4**

(i) $K_{22} \in \mathcal{L}^+(V)$, $K_{12} \in \mathcal{L}(U, V)$ is injective and $K_{21} = K_{12}^*$.
(ii) There is a strictly increasing sequence $(m_k)_{k \in \mathbb{N}} \subset \mathbb{N}$ such that for each $n \in \mathbb{N}$ with $n \leq m_k$, it holds

$$K_{22}(B_V^n) \subset B_V^{m_k}, \quad K_{12}(B_U^n) \subset B_V^{m_k} \quad \text{and} \quad K_{21}(B_V^n) \subset B_U^{m_k}.$$

(iii) There is a constant $c_K \in (0, \infty)$ such that

$$(K_{12} K_{21} v, v)_V \leq c_K (K_{22} v, v)_V.$$

(iv) $\|D\Phi\|_{L^\infty(\mu_1)}^2 < \infty$.

**Remark 3.5** Let $n \in \mathbb{N}$ and $f = \varphi(P_n(\cdot), P_n(\cdot)) \in \mathcal{F}C_b^\infty(B_W)$ be given. One can show that $D_2 f = \sum_{i=1}^n \partial_{i,2} \varphi(P_n(\cdot), P_n(\cdot)) e_i \in B_V^n$ and $D_1 f \in B_U^n$. Hence the invariance properties stated in Assumption 3.4 ensure that expressions like $Q_2^{-1} K_{22} D_2 f$, $Q_1^{-1} K_{21} D f$ are well-defined.

**Theorem 3.6** [16, Lemma 4.5 and Theorem 4.11] *Suppose Assumption 3.4 is valid. Then $(S, \mathcal{F}C_b^\infty(B_W))$ is symmetric, $(A_\Phi, \mathcal{F}C_b^\infty(B_W))$ is antisymmetric and $(L_\Phi, \mathcal{F}C_b^\infty(B_W))$ is essentially m-dissipative in $L^2(\mu^\Phi)$. Moreover we have for all $f \in \mathcal{F}C_b^\infty(B_W)$*

$$-(Sf, f)_{L^2(\mu^\Phi)} = \int_W (K_{22} D_2 f, D_2 f)_V \, d\mu^\Phi \quad \text{and} \quad \int_W L_\Phi f \, d\mu^\Phi = 0.$$

*In particular, the measure $\mu^\Phi$ is invariant for $(L_\Phi, \mathcal{F}C_b^\infty(B_W))$.*

**Remark 3.7** Since symmetric, antisymmetric and dissipative operators are closable, it is reasonable to denote by $(S, D(S))$, $(A_\Phi, D(A_\Phi))$ and $(L_\Phi, D(L_\Phi))$ the closures of the operators $(S, \mathcal{F}C_b^\infty(B_W))$, $(A_\Phi, \mathcal{F}C_b^\infty(B_W))$ and $(L_\Phi, \mathcal{F}C_b^\infty(B_W))$, respectively.

In order to construct a weak solution in form of a diffusion process with infinite lifetime we have to verify the assumption below. The construction is done in Theorem 4.3 and is based on the existence of a compact $\mu^\Phi$-nest.

**Assumption 3.8**

(i) *The operator $K_{22}$ is of trace class and $Q_2^{-1}K_{22}$ restricted to span$\{e_1, ..., e_n\}$ is non-negative for all $n \in \mathbb{N}$.*
(ii) *There exists $\rho \in L^1(\mu^\Phi)$, such that for all $n \in \mathbb{N}$, the function $\rho_n$ defined by*

$$W \ni (u, v) \to \rho_n(u, v) = \rho(\overline{P}_n(u), \overline{P}_n(v)) \in \mathbb{R},$$

*is in $L^1(\mu^\Phi)$. Furthermore for all $(u, v) \in W$ we have*

$$(\overline{P}_n(u), Q_1^{-1}K_{21}\overline{P}_n(v))_U - (\overline{P}_n(v), Q_2^{-1}K_{12}\overline{P}_n(u))_V \leq \rho(\overline{P}_n(u), \overline{P}_n(v)),$$

*and $(\rho_n)_{n\in\mathbb{N}}$ converges to $\rho$ in $L^1(\mu^\Phi)$.*

Below we introduce an infinite-dimensional Ornstein–Uhlenbeck operator perturbed by the gradient of the potential $\Phi$. This is necessary, since we plan to use the general abstract hypocoercivity framework, described by Grothaus and Stilgenbauer. Indeed, in Sect. 6 we see that these operators naturally appear as we show exponential convergence to equilibrium of the semigroup generated by $(L_\Phi, D(L_\Phi))$.

**Definition 3.9** Define $(C, D(C))$ and $(Q_1^{-1}C, D(Q_1^{-1}C))$ on $U$ by

$$C = K_{21}Q_2^{-1}K_{12} \quad \text{with} \quad D(C) = \{u \in U \mid K_{12}u \in D(Q_2^{-1})\} \quad \text{and}$$
$$Q_1^{-1}C = Q_1^{-1}K_{21}Q_2^{-1}K_{12} \quad \text{with} \quad D(Q_1^{-1}C) = \{u \in D(C) \mid Cu \in D(Q_1^{-1})\}.$$

Moreover we define the operator $(N, \mathcal{F}C_b^\infty(B_U))$ by

$$\mathcal{F}C_b^\infty(B_U) \ni f \mapsto Nf = \text{tr}[CD^2 f] - (\cdot, Q_1^{-1}CDf)_U - (D\Phi, CDf)_U$$
$$\in L^2(U, \mu_1^\Phi, \mathbb{R}).$$

We want to highlight that the definition of the infinite-dimensional Ornstein Uhlenbeck operator above, allows potentially unbounded diffusion operators $(C, D(C))$ as coefficients. One can check that $(C, D(C))$ is symmetric and positive. Moreover, for each $n \in \mathbb{N}$ with $n \leq m_k$, it holds $B_U^n \subset D(C)$ as well as $C(B_U^n) \subset B_U^{m_k}$.
Assumption 3.10 below is central to apply the abstract hypocoercivity concept described in Sect. 5 and therefore to obtain exponential convergence to equilibrium of the semigroup generated by $(L_\Phi, D(L_\Phi))$. In particular, the assumption guarantees essential m-dissipativity of the infinite-dimensional Ornstein–Uhlenbeck operators $(N, \mathcal{F}C_b^\infty(B_U))$ and corresponding regularity estimates, compare Theorem 3.11.

**Assumption 3.10**

(i) $\|D\Phi\|_{L^\infty(\mu_1)}^2 < \frac{1}{4\lambda_1}$.
(ii) $\Phi : U \to (-\infty, \infty]$ *is convex and lower semicontinuous.*

(iii) *There exists a constant $\omega_1 \in (0, \infty)$ s.t. for all $f \in \mathcal{F}C_b^\infty(B_V)$*

$$\int_V (K_{22}Df, Df)_V \, d\mu_2 \geq \omega_1 \int_V \left(f - (f, 1)_{L^2(\mu_2)}\right)^2 d\mu_2.$$

(iv) *There is a constant $\omega_2 \in (0, \infty)$ s.t. for all $f \in \mathcal{F}C_b^\infty(B_U)$*

$$\int_U (Q_2^{-1} K_{12} Df, K_{12} Df)_V \, d\mu_1^\Phi \geq \omega_2 \int_U \left(f - (f, 1)_{L^2(\mu_1^\Phi)}\right)^2 d\mu_1^\Phi.$$

(v) *Assume that there is a constant $C_1 \in (0, \infty)$ s.t. for all $f \in \mathcal{F}C_b^\infty(B_W)$ and $g = (Id - PA_\Phi^2 P)f$, where $P$ is defined in Definition 6.1, it holds*

$$\int_U (K_{22} Q_2^{-1} K_{12} D_1 f_S, Q_2^{-1} K_{22} Q_2^{-1} K_{12} D_1 f_S)_V \, d\mu_1^\Phi \leq C_1^2 \|g\|_{L^2(\mu^\Phi)}^2.$$

**Theorem 3.11** [16, Theorem 3.6 and 3.11] *Assume Item (i) from Assumption 3.10. Then the operator $(N, \mathcal{F}C_b^\infty(B_U))$ is essentially m-dissipative in $L^2(U, \mu_1^\Phi, \mathbb{R})$. Moreover, we have*

$$(Nf, g)_{L^2(U, \mu_1^\Phi, \mathbb{R})} = \int_U -(CDf, Dg)_U \, d\mu_1^\Phi \quad \text{for all} \quad f, g \in \mathcal{F}C_b^\infty(B_U).$$

*We denote the closure of $(N, \mathcal{F}C_b^\infty(B_U))$ by $(N, D(N))$ and the corresponding resolvent in $\alpha \in (0, \infty)$ by $R(\alpha, N)$. Suppose in addition that Item (ii) is valid. Then it holds for $g = \alpha f - Nf$, where $\alpha \in (0, \infty)$ and $f \in \mathcal{F}C_b^\infty(B_U)$*

$$\int_U \alpha f^2 + (CDf, Df)_U \, d\mu_1^\Phi \leq \frac{1}{\alpha} \int_U g^2 \, d\mu_1^\Phi \quad \text{and}$$

$$\int_U \operatorname{tr}[(CD^2 f)^2] + \|Q_1^{-\frac{1}{2}} CDf\|_U^2 \, d\mu_1^\Phi \leq 4 \int_U g^2 \, d\mu_1^\Phi.$$

In applications, compare Sect. 7, the regularity estimates above can be used to check Item (v) from Assumption 3.10. In order to check the Poincaré type inequalities stated in Assumption 3.10 above the general Poincaré inequality from the proposition below is useful.

**Proposition 3.12** *Suppose $\Phi : U \to (-\infty, \infty]$ is convex, bounded from below, lower semicontinuous and not identically to $\infty$. Then for all $f \in \mathcal{F}C_b^\infty(B_U)$ it holds*

$$\int_U (Q_1 Df, Df)_U \, d\mu_1^\Phi \geq \lambda_1 \int_U \left(f - (f, 1)_{L^2(\mu_1^\Phi)}\right)^2 d\mu_1^\Phi.$$

**Proof** The idea of the proof is to approximate $\Phi$ by a sequence of convex and smooth functions. Afterwards, we apply the Poincaré inequality from [1, Proposition 4.5].

I.e., let $\alpha > 0$ be given and denote by $\Phi_\alpha$ the Moreau-Yosida approximation, defined by

$$\Phi_\alpha(u) = \inf \left\{ \Phi(\bar{u}) + \frac{1}{2\alpha} \|u - \bar{u}\|_U^2 \mid \bar{u} \in U \right\}, \quad u \in U.$$

It is well known, see [12], that $\Phi_\alpha$ is convex, differentiable with Lipschitz continuous derivative. Furthermore, $-\infty < \inf_{\bar{u} \in U} \Phi(\bar{u}) \leq \Phi_\alpha(u) \leq \Phi(u)$, for all $u \in U$ and $\lim_{\alpha \to 0} \Phi_\alpha(u) = \Phi(u)$, for all $u \in U$. But to use the Poincaré inequality from [1, Proposition 4.5] $\Phi_\alpha$ is not regular enough. Therefore we take $\beta > 0$ and define the function $\Phi_{\alpha,\beta}$ by

$$\Phi_{\alpha,\beta}(u) = \int_U \Phi_\alpha(e^{\beta B} u + \bar{u}) d\mu_\beta(\bar{u}), \quad u \in U,$$

where $(B, D(B))$ is a self-adjoint negative definite operator such that $B^{-1}$ is of trace class and $\mu_\beta$ the infinite dimensional centered non-degenerate Gaussian measure with covariance operator $\frac{1}{2} B^{-1}(e^{2\beta B} - Id)$. Note that $\Phi_{\alpha,\beta}$ is defined in terms of the Mehler formula of an infinite-dimensional Ornstein–Uhlenbeck operator, see [8, Section 8.3]. Recalling the smoothing properties of such Ornstein–Uhlenbeck semigroups, one can check, compare [8, Section 11.6], that $\Phi_{\alpha,\beta}$ is convex and has derivatives of all orders. Moreover, $D\Phi_{\alpha,\beta}$ is Lipschitz continuous and has bounded derivatives of all orders. This is enough to verify that $\Phi_{\alpha,\beta}$ fulfills Hypothesis 1.3 from [1]. Hence, by the Poincaré inequality from [1, Proposition 4.5] we obtain for all $f \in \mathcal{F}C_b^\infty(B_U)$

$$\int_U (Q_1 Df, Df)_U d\mu_1^{\Phi_{\alpha,\beta}} \geq \lambda_1 \int_U \left( f - (f, 1)_{L^2(\mu_1^{\Phi_{\alpha,\beta}})} \right)^2 d\mu_1^{\Phi_{\alpha,\beta}}. \qquad (2)$$

Since the derivative of $\Phi_\alpha$ is Lipschitz continuous one can show that $\Phi_\alpha$ has at most quadratic growth. Hence, we obtain $\lim_{\beta \to 0} \Phi_{\alpha,\beta}(u) = \Phi_\alpha(u)$, for all $u \in U$, by the theorem of dominated convergence. In particular, we get $\lim_{\alpha \to 0} \lim_{\beta \to 0} \Phi_{\alpha,\beta}(u) = \Phi(u)$, for all $u \in U$. An application of the dominated convergence theorem results in

$$\lim_{\alpha \to 0} \lim_{\beta \to 0} \mu_1^{\Phi_{\alpha,\beta}} = \mu_1^\Phi \quad \text{weakly.}$$

I.e., taking the limits $\beta \to 0$ and $\alpha \to 0$ in Inequality (2) yields the claim. □

Note that the Poincaré inequality from the proposition above is valid without assuming that $\Phi \in W^{1,2}(U, \mu_1, \mathbb{R})$ and generalizes the one from [11, Theorem 12.3.8]. Indeed, since $Q_1$ is trace class, it is easier to estimate $\lambda_1 \int_U \left( f - (f, 1)_{L^2(\mu_1^\Phi)} \right)^2 d\mu_1^\Phi$ by $\int_U (Df, Df)_U d\mu_1^\Phi$ from above, as it is done in [11], than by $\int_U (Q_1 Df, Df)_U d\mu_1^\Phi$. Nevertheless, the proof of Proposition 3.12 is inspired by the one from [11, Theorem 12.3.8].

## 4 Weak solutions

In this section we construct a solution to the infinite-dimensional Langevin equation

$$dU_t = K_{21} Q_2^{-1} V_t dt \tag{3}$$
$$dV_t = -(K_{22} Q_2^{-1} V_t + K_{12} Q_1^{-1} U_t + K_{12} D\Phi(U_t)) dt + \sqrt{2K_{22}} dW_t, \tag{4}$$

where the involved operators were specified in Assumption 3.4, above. As mentioned in the introduction, this equation corresponds to the infinite-dimensional Langevin operator $(L_\Phi, \mathcal{F}C_b^\infty(B_W))$. First we construct a martingale solution, compare Proposition 4.4 below. Second we show that the martingale solution also provides a weak solution.

Let $\mathcal{T}$ be the weak topology on $(W, (\cdot, \cdot)_W)$. Clearly $(W, \mathcal{T})$ is a Lusin space (i.e. the continuous one-to-one image of a Polish space). Denote by $\mathcal{B}_b(W)$, $C([0, \infty), W)$ and $C_b(W)$ the space of bounded $\mathcal{B}(W)$ measurable functions from $W$ to $\mathbb{R}$, the space of continuous functions from $[0, \infty)$ to $W$ and the space of continuous bounded functions from $W$ to $\mathbb{R}$, respectively. The potential theoretic notions we need in this article are from [4] and [23]. During the whole section we assume Assumption 3.4. Therefore, we know by Theorem 3.6 we know that the closure $(L_\Phi, D(L_\Phi))$ of $(L_\Phi, \mathcal{F}C_b^\infty(B_W))$ in $L^2(\mu^\Phi)$ generates a strongly continuous contraction semigroup $(T_t)_{t \geq 0}$.

**Remark 4.1** Using that $\mu^\Phi$ is invariant for $L_\Phi$ and [15, Lemma 1.9, App. B] we know that the semigroup $(T_t)_{t \geq 0}$ generated by $(L_\Phi, D(L_\Phi))$ is Markovian, i.e. positive preserving and $T_t 1 = 1$ for all $t \geq 0$. Note that the same holds true for the strongly continuous contraction resolvent $(R(\alpha, L_\Phi))_{\alpha > 0}$ corresponding to $(L_\Phi, D(L_\Phi))$.

For $k, n \in \mathbb{N}$ and $(u, v) \in W$ we set:

$$N(u, v) = \|(u, v)\|_W^2 = \|u\|_U^2 + \|v\|_V^2,$$
$$N_n(u, v) = N(\overline{P}_n(u), \overline{P}_n(v)) \quad \text{and}$$
$$F_k = \{(u, v) \in W \mid N(u, v) \leq k\}.$$

Obviously, for all $k \in \mathbb{N}$ the sets $F_k$ are $\mathcal{T}$-compact and increasing. Below we show that the sequence $(F_k)_{k \in \mathbb{N}}$ provides a $\mathcal{T}$-compact $\mu^\Phi$-nest, if Assumption 4.2 is valid. The idea of the proof follows the lines of [4, Proposition 5.5].

**Lemma 4.2** *Assume Assumption 3.8 holds true. Then the sequence of functions $(h_n)_{n \in \mathbb{N}}$ defined by*

$$h_n(u, v) = (\overline{P}_n(v), Q_2^{-1} K_{22} \overline{P}_n(v))_V, \quad (u, v) \in W,$$

*is non-negative and converges to a non-negative function $h$ in $L^1(\mu^\Phi)$. Moreover, the increasing sequence of $\mathcal{T}$-compact sets $(F_k)_{k \in \mathbb{N}} = (\{(u, v) \in W \mid N(u, v) \leq k\})_{k \in \mathbb{N}}$ is a $\mu^\Phi$-nest.*

**Proof** By Assumption 3.8 Item (*i*) we know that $(h_n)_{n\in\mathbb{N}}$ is a sequence of non-negative functions. Moreover, we can calculate for $m, n \in \mathbb{N}$ with $m \geq n$

$$\int_W |h_m(u,v) - h_n(u,v)| \mathrm{d}\mu^\Phi(u,v)$$

$$= \int_W \left| \sum_{i,j=n+1}^m (e_i, Q_2^{-1} K_{22} e_j)_V (v, e_i)_V (v, e_j)_V \right| \mathrm{d}\mu^\Phi(u,v)$$

$$= \sum_{i,j=n+1}^m (e_i, Q_2^{-1} K_{22} e_j)_V \int_W (v, e_i)_V (v, e_j)_V \mathrm{d}\mu^\Phi(u,v)$$

$$= \sum_{i=n+1}^m (e_i, K_{22} e_i)_V.$$

Note that we used Lemma 3.3 in the calculation above. Since we assume that $K_{22}$ is of trace class, we observe that $(h_n)_{n\in\mathbb{N}}$ is a non-negative Cauchy-sequence in $L^1(\mu^\Phi)$ and therefore convergent to a non-negative function $h \in L^1(\mu^\Phi)$. Using a suitable selection of cut-off functions we obtain $N_n \in D(L_\Phi)$ for all $n \in \mathbb{N}$ with

$$L_\Phi N_n(u,v)$$
$$= 2\sum_{i=1}^n (K_{22} e_i, e_i)_V - 2(\overline{P}_n(v), Q_2^{-1} K_{22} \overline{P}_n(v))_V - 2(\overline{P}_n(u), Q_1^{-1} K_{21} \overline{P}_n(v))_U$$
$$+ 2(\overline{P}_n(v), Q_2^{-1} K_{12} \overline{P}_n(u))_V - 2(D\Phi(u), K_{21} \overline{P}_n(v))_U$$
$$\geq -2h_n(v) - 2\rho_n(u,v) - 2\|K_{12} D\Phi\|_{L^\infty(\mu_1^\Phi)} \|\overline{P}_n(v)\|_V.$$

Setting $g(u,v) = N(u,v) + 2\rho(u,v) + 2\|K_{12} D\Phi\|_{L^\infty(\mu_1^\Phi)} \|v\|_V$ and $g_n(u,v) = g(\overline{P}_n(u), \overline{P}_n(v))$ we get for all $n \in \mathbb{N}$

$$(Id - L_\Phi) N_n \leq g_n + 2h_n.$$

Similar to [4, Proposition 5.5] one can show that the closure of $(L_\Phi, \mathcal{F}C_b^\infty(B_W))$ in $L^1(\mu^\Phi)$ also generates a strongly continuous contraction semigroup with corresponding resolvent denoted by $(R_1(\alpha, L_\Phi))_{\alpha>0}$. Note that this resolvent coincides with $(R(\alpha, L_\Phi))_{\alpha>0}$ on $L^2(\mu^\Phi)$. Applying $R_1(1, L_\Phi)$ on both sides of the inequality above results in

$$N_n \leq R_1(1, L_\Phi) g_n + 2R_1(1, L_\Phi) h_n.$$

Note that the inequality should be read $\mu^\Phi$-a.e. and makes use of the Markovianity of the resolvent. As $(N_n)_{n\in\mathbb{N}}$, $(g_n)_{n\in\mathbb{N}}$ and $(h_n)_{n\in\mathbb{N}}$ converge to $N$, $g$ and $h$ in $L^1(\mu^\Phi)$ we obtain $\mu^\Phi$-a.e.

$$N(u,v) \leq R_1(1, L_\Phi) g + 2R_1(1, L_\Phi) h.$$

Proceeding as in [4, Proposition 5.5] shows the assertion. □

Having the $\mathcal{T}$-compact $\mu^\Phi$-nest at hand the application of [4, Theorem 1.1] is possible.

**Theorem 4.3** *Assume that Assumption 3.8 holds true. Then there exists a countable $\mathbb{Q}$-algebra $\mathcal{A} \subset D(L_\Phi) \cap C_b(W)$ such that*

1. *$\mathcal{A}$ is a core for $(L_\Phi, D(L_\Phi))$ and*
2. *$\mathcal{A}$ separates the points of $W$.*

*Let $\mathcal{T}_0$ be the topology on $W$ generated by $\mathcal{A}$. Then there exists a $\mu^\Phi$-standard right process*

$$X = (\Omega, \mathcal{F}, (\mathcal{F}_t)_{t\geq 0}, (X_t)_{t\geq 0}, (\theta_t)_{t\geq 0}, (\mathbb{P}^w)_{w\in W}),$$

*(see [4, Appendix B.]) with the state space $W$ (endowed with the topology $\mathcal{T}_0$) whose transition semigroup denoted by $(p_t)_{t>0}$ coincides with $(T_t)_{t>0}$ on $L^2(\mu^\Phi)$. The transition semigroup $(p_t)_{t>0}$ is defined via*

$$p_t g(w) = \int_\Omega g(X_t(\omega))\mathrm{d}\mathbb{P}^w(\omega), \quad w \in W \text{ and } g \in B_b(W),$$

*and identified with its extension to $L^2(\mu^\Phi)$.*

*Denote by $\mathbb{P}_{\mu^\Phi}$ the probability measure on $(\Omega, \mathcal{F})$ defined for $A \in \mathcal{F}$ by*

$$\mathbb{P}_{\mu^\Phi}(A) = \int_W \mathbb{P}^w(A)\mathrm{d}\mu^\Phi(w), \quad A \in \mathcal{F}.$$

*Then it holds*

(i) *The process has infinite lifetime $\mathbb{P}_{\mu^\Phi}$-a.e..*
(ii) *The process is $\mathcal{T}$-continuous, $\mathbb{P}_{\mu^\Phi}$-a.e..*
(iii) *Every element from $D(L_\Phi)$ has a $\mu^\Phi$-quasi continuous version (with respect to the topology $\mathcal{T}_0$).*

*By (i) and (ii), $X$ is a diffusion process with infinite lifetime.*

**Proof** Note that $\bigcup_{n\in\mathbb{N}} \mathcal{F}C_b^\infty(B_W, n) = \mathcal{F}C_b^\infty(B_W) \subset D(L_\Phi) \cap C_b(W)$ is a core for $(L_\Phi, D(L_\Phi))$. Now given $f = \varphi(P_n(\cdot), P_n(\cdot)) \in \mathcal{F}C_b^\infty(B_W, n)$ for some $n \in \mathbb{N}$. Using a suitable sequence of cut-off functions provided by Lemma 3.2 we can approximate $f$ by a function $g = \psi(P_n(\cdot), P_n(\cdot))$, where $\psi \in C_0^\infty(\mathbb{R}^n \times \mathbb{R}^n)$, in $(L_\Phi, D(L_\Phi))$ graph-norm. As there exists a countable set $\mathcal{A}_n = \{\psi_m^n \in C_0^\infty(\mathbb{R}^n \times \mathbb{R}^n) \mid m \in \mathbb{N}\}$ which is dense in $C_0^\infty(\mathbb{R}^n \times \mathbb{R}^n)$ w.r.t. the norm $|\cdot|_\infty + |D\cdot|_\infty + |D^2\cdot|_\infty$, we find a sequence $(\psi_m^n)_{m\in\mathbb{N}} \subset \mathcal{A}_n$ converging to $\psi$ w.r.t. $|\cdot|_\infty + |D\cdot|_\infty + |D^2\cdot|_\infty$. We obtain convergence of the sequence $(\psi_m^n(P_n(\cdot), P_n(\cdot)))_{m\in\mathbb{N}}$ to $g$ w.r.t. the $(L_\Phi, D(L_\Phi))$ graph-norm. We can conclude that $\mathcal{A}_n$ is dense in $\mathcal{F}C_b^\infty(B_W, n)$ w.r.t. the $(L_\Phi, D(L_\Phi))$ graph-norm. In particular, the smallest $\mathbb{Q}$-algebra $\mathcal{A}$ containing

$$\bigcup_{n\in\mathbb{N}} \{\varphi(P_n(\cdot), P_n(\cdot)) \mid \varphi \in \mathcal{A}_n\},$$

is a core for $(L_\Phi, D(L_\Phi))$. Moreover, $\mathcal{A}$ is a subset of $D(L_\Phi) \cap C_b(W)$ and by construction countable. It is easy to see that $\mathcal{A}$ also separates the points of $W$.

The existence of an $\mathcal{T}$-compact $\mu^\Phi$-nest is shown in Lemma 4.2. By [4, Theorem 1.1], [7, Lemma 2.2.8] and [23, Exercise IV.2.9] there is a $\mu^\Phi$-standard right process with state space $W$ (endowed with the topology $\mathcal{T}_0$) whose transition semigroup $(p_t)_{t>0}$ extended to $L^2(\mu^\Phi)$ coincides with $(T_t)_{t>0}$. Additionally by [4, Theorem 1.1] we know that Item $(iii)$ is fulfilled and that $(X_t)_{t\geq 0}$ is càdlàg w.r.t. the weak topology $\mathcal{T}$, $\mathbb{P}_{\mu^\Phi}$-a.e.. Adapting the arguments from [4, Proposition 5.6.] we obtain $\mathcal{T}$-continuity of $(X_t)_{t\geq 0}$.

Infinite lifetime $\mathbb{P}_{\mu^\Phi}$-a.e. follows as in [7, Theorem 6.4.2.], since $\mu^\Phi$ is invariant for $(L_\Phi, \mathcal{F}C_b^\infty(B_W))$, by Theorem 3.6. □

Denote by $(U_t)_{t\geq 0}$ and $(V_t)_{t\geq 0}$ the projection of $(X_t)_{t\geq 0}$ to $U$ and $V$, respectively. Invoking [4, Proposition 1.4] and [7, Lemma 2.1.8] the following proposition is immediate.

**Proposition 4.4** [4, Proposition 1.4] and [7, Lemma 2.1.8] *Suppose we are in the situation of Theorem* 4.3 *and*

$$X = (\Omega, \mathcal{F}, (\mathcal{F}_t)_{t\geq 0}, (X_t)_{t\geq 0}, (\theta_t)_{t\geq 0}, (\mathbb{P}^w)_{w\in W}),$$

*is the $\mu^\Phi$-standard right process provided there. Given a non-negative $g_0 \in L^2(\mu^\Phi)$ with $\int_W g_0 \mathrm{d}\mu^\Phi = 1$. Denote by $\nu^\Phi$ the measure with $\mu^\Phi$-density $g_0$ and define the probability measure $\mathbb{P}_{\nu^\Phi}$ on $(\Omega, \mathcal{F})$ by*

$$\mathbb{P}_{\nu^\Phi}(A) = \int_W \mathbb{P}^w(A) \mathrm{d}\nu^\Phi(w), \quad A \in \mathcal{F}.$$

*Then X solves the martingale problem for $(L_\Phi, D(L_\Phi))$ under $\mathbb{P}_{\nu^\Phi}$. I.e., for all $f \in D(L_\Phi)$ the process $(M_t^{[f]})_{t\geq 0}$ defined by*

$$M_t^{[f]} = f(U_t, V_t) - f(U_0, V_0) - \int_{[0,t]} L_\Phi f(U_s, V_s) \mathrm{d}\lambda(s), \quad t \geq 0, \qquad (5)$$

*defines an $(\mathcal{F}_t)_{t\geq 0}$ martingale w.r.t. $\mathbb{P}_{\nu^\Phi}$.*

*If in addition $f^2 \in D(L_\Phi)$ and $L_\Phi f \in L^4(\mu^\Phi)$ then the process $(N_t^{[f]})_{t\geq 0}$ defined by*

$$N_t^{[f]} = (M_t^{[f]})^2 - \int_{[0,t]} (L_\Phi f^2)(U_s, V_s) - 2(f L_\Phi f)(U_s, V_s) \mathrm{d}\lambda(s), \quad t \geq 0, \quad (6)$$

*also defines an $(\mathcal{F}_t)_{t\geq 0}$ martingale w.r.t. $\mathbb{P}_{\nu^\Phi}$.*

**Remark 4.5** Given $f \in D(L_\Phi)$ and $t \in [0, \infty)$. In [6, Lemma 5.1] it is shown, that the random variables $M_t^{[f]}$ and $N_t^{[f]}$ in the proposition above are well-defined, i.e. $\mathbb{P}_{\nu^\Phi}$-a.s. independent of the representative of $f$ and $L_\Phi f$ in $L^2(\mu^\Phi)$.

To construct weak solutions we calculate how the operator $(L_\Phi, D(L_\Phi))$ acts on a certain class of functions.

**Lemma 4.6** *For $i, j \in \mathbb{N}$ we define*

$$W \ni (u, v) \mapsto f_i(u, v) = (u, d_i)_U \in \mathbb{R} \quad \text{and}$$
$$W \ni (u, v) \mapsto g_i(u, v) = (v, e_i)_V \in \mathbb{R}.$$

*Then it holds $f_i, g_i, f_i f_j, g_i g_j \in D(L_\Phi)$, $L_\Phi f_i, L_\Phi g_i \in L^4(\mu^\Phi)$ and*

$$L_\Phi f_i = (v, Q_2^{-1} K_{12} d_i)_V$$
$$L_\Phi f_i f_j = f_j L_\Phi f_i + f_i L_\Phi f_j$$
$$L_\Phi g_i = -(v, Q_2^{-1} K_{22} e_i)_V - (u, Q_1^{-1} K_{21} e_i)_U - (D\Phi(u), K_{21} e_i)_U$$
$$L_\Phi g_i g_j = 2(e_j, K_{22} e_i)_V + g_j L_\Phi g_i + g_i L_\Phi g_j.$$

*Proof* Use suitable cut-off functions provided by Lemma 3.2 and note that $D\Phi$ is bounded by Assumption 3.4. □

Using the process $X = (\Omega, \mathcal{F}, (\mathcal{F}_t)_{t \geq 0}, (X_t)_{t \geq 0}, (\theta_t)_{t \geq 0}, (\mathbb{P}^w)_{w \in W})$ provided by Theorem 4.3 we can construct a weak solution to the infinite-dimensional Langevin equation.

**Theorem 4.7** *There exists a cylindrical Brownian motion $(W_t)_{t \geq 0}$ on $(\Omega, \mathcal{F}, \mathbb{P}_{\mu^\Phi})$ with values in $V$ s.t. it holds $\mathbb{P}_{\mu^\Phi}$-a.s. for all $\vartheta \in D(Q_2^{-1} K_{12})$, $\theta \in D(Q_1^{-1} K_{21}) \cap D(Q_2^{-1} K_{22})$ and $t \in [0, \infty)$*

$$(U_t, \vartheta)_U = (U_0, \vartheta)_U + \int_{[0,t]} (V_s, Q_2^{-1} K_{12} \vartheta)_V d\lambda(s) \quad \text{and}$$
$$(V_t, \theta)_V = (V_0, \theta)_V - \int_{[0,t]} (V_s, Q_2^{-1} K_{22} \theta)_V + (U_s, Q_1^{-1} K_{21} \theta)_U$$
$$+ (D\Phi(U_s), K_{21} \theta)_U d\lambda(s) + (\sqrt{2 K_{22}} W_t, \theta)_V.$$

*Proof* Given $i, j \in \mathbb{N}$ and $t \in [0, \infty)$. By the calculations in Lemma 4.6 and Equation (6) we obtain

$$N_t^{[g_i]} = (M_t^{[g_i]})^2 - 2(K_{22} e_i, e_i)_V t \quad \text{and}$$
$$N_t^{[g_i + g_j]} = (M_t^{[g_i + g_j]})^2 - 4(K_{22} e_i, e_j)_V t - 2(K_{22} e_i, e_i)_V t - 2(K_{22} e_j, e_j)_V t.$$

Hence the quadratic variation of $M^{[g_i]}$ and $M^{[g_j]}$ at time $t$ is given by

$$\langle M^{[g_i]}, M^{[g_j]} \rangle_t = 2(K_{22} e_i, e_j)_V t.$$

Now define a $V$ valued process $M^* = (M_t^*)_{t \geq 0}$ by

$$M_t^* = \sum_{i=1}^{\infty} M_t^{[g_i]} e_i, \quad t \in [0, \infty).$$

Using that $K_{22}$ is trace class we can argue as in [26, Proposition 2.1.10] to show that the series above converges in $L^2(\Omega, \mathcal{F}, \mathbb{P}_{\nu^\Phi}, C([0, \infty), W))$ and defines a $2K_{22}$ Brownian motion. Furthermore, one can show the existence of a cylindrical Brownian motion $(W_t)_{t \geq 0}$ such that $M_t^* = \sqrt{2K_{22}} W_t$, for all $t \in [0, \infty)$. Hence, we obtain by the calculations in Lemma 4.6 and Equation (5), $\mathbb{P}_{\mu^\Phi}$-a.s. for $t \in [0, \infty)$

$$(\sqrt{2K_{22}} W_t, e_i)_V = M_t^{[g_i]} = (V_t, e_i)_V - (V_0, e_i)_V$$
$$+ \int_{[0,t]} (V_s, Q_2^{-1} K_{22} e_i)_V + (U_s, Q_1^{-1} K_{21} e_i)_U + (D\Phi(U_s), K_{21} e_i)_U \, d\lambda(s).$$

Now given an element $\theta \in D(Q_1^{-1} K_{21}) \cap D(Q_2^{-1} K_{22})$. If we multiply the equation above with $(\theta, e_i)_V$ and sum over all $i \in \mathbb{N}$ we get $\mathbb{P}_{\mu^\Phi}$-a.s.

$$(\sqrt{K_{22}} W_t, \theta)_V$$
$$= (V_t, \theta)_V - (V_0, \theta)_V + \lim_{i \to \infty} \int_{[0,t]} (V_s, Q_2^{-1} K_{22} \overline{P}_{m_i}(\theta))_V$$
$$+ (U_s, Q_1^{-1} K_{21} \overline{P}_{m_i}(\theta))_U + (D\Phi(U_s), K_{21} \overline{P}_{m_i}(\theta))_U \, d\lambda(s)$$
$$= (V_t, \theta)_V - (V_0, \theta)_V + \int_{[0,t]} (V_s, Q_2^{-1} K_{22} \theta)_V + (U_s, Q_1^{-1} K_{21} \theta)_U$$
$$+ (D\Phi(U_s), K_{21} \theta)_U \, d\lambda(s).$$

In the last equation we used the dominated convergence theorem. It is applicable since for $\theta \in D(Q_1^{-1} K_{21}) \cap D(Q_2^{-1} K_{22})$ we have pointwisely for $i \to \infty$

$$(V_s, Q_2^{-1} K_{22} \overline{P}_{m_i}(\theta))_V = (V_s, \overline{P}_{m_i}(Q_2^{-1} K_{22} \theta))_V \to (V_s, Q_2^{-1} K_{22} \theta)_V,$$
$$(U_s, Q_1^{-1} K_{21} \overline{P}_{m_i}(\theta))_U = (U_s, \overline{P}_{m_i}(Q_1^{-1} K_{21} \theta))_U \to (U_s, Q_1^{-1} K_{21} \theta)_U \text{ and}$$
$$(D\Phi(U_s), K_{21} \overline{P}_{m_i}(\theta))_U \to (D\Phi(U_s), K_{21} \theta)_U.$$

Moreover, $(V_s)_{s \in [0,t]}$ and $(U_s)_{s \in [0,t]}$ are bounded $\mathbb{P}_{\mu^\Phi}$-a.s., since they are $\mathcal{T}$-continuous $\mathbb{P}_{\mu^\Phi}$-a.s.. The quadratic variation of $M^{[f_i]}$ at time $t \in [0, \infty)$ is equal to zero, by the calculations in Lemma 4.6 and Equation (6). Therefore we get $\mathbb{P}_{\mu^\Phi}$-a.s. for $t \in [0, \infty)$

$$0 = M_t^{[f_i]} = (U_t, d_i)_U - (U_0, d_i)_U - \int_{[0,t]} (V_s, Q_2^{-1} K_{12} d_i)_V \, d\lambda(s).$$

Using similar arguments as above we obtain $\mathbb{P}_{\mu^\Phi}$-a.s. for $\vartheta \in D(Q_2^{-1}K_{12})$ and for $t \in [0, \infty)$

$$(U_t, \vartheta)_U = (U_0, \vartheta)_U + \int_{[0,t]} (V_s, Q_2^{-1}K_{12}\vartheta)_V \, d\lambda(s).$$

□

## 5 The Hilbert space hypocoercivity method

This section is devoted to the abstract Hilbert space hypocoercivity method presented in [17]. It is a fundamental tool to establish hypocoercivity of infinite-dimensional Langevin dynamics. The method originated by J. Dolbeault, C. Mouhot and C. Schmeiser in [13], where an algebraic hypocoercivity method, i.e. excluding domain issues, for linear kinetic equations is developed.

Below, $W$ always denotes a real Hilbert space with scalar product $(\cdot, \cdot)$ and induced norm $\|\cdot\|$. All considered operators are assumed to be linear and defined on linear subspaces of $W$.

**Hypocoercivity data (D).**

(D1) *The Hilbert space.* Let $(W, \mathcal{F}, \mu)$ be a probability space. Set $H = L^2(W, \mu, \mathbb{R})$ and equip it with the usual standard scalar product $(\cdot, \cdot)$ and denote by $\|\cdot\|$ the induced norm.

(D2) *The $C_0$-semigroup and its generator $L$.* $(L, D(L))$ is a linear operator on $H$ generating a strongly continuous semigroup $(T_t)_{t \geq 0}$.

(D3) *Core property of $L$.* Let $D \subset D(L)$ be a dense subspace of $H$ which is an operator core for $(L, D(L))$.

(D4) *Decomposition of $L$.* Let $(S, D(S))$ be symmetric and let $(A, D(A))$ be closed and antisymmetric on $H$ such that $D \subset D(S) \cap D(A)$ as well as $L_{|D} = S - A$.

(D5) *Orthogonal projection.* Let $P : H \to H$ be an orthogonal projection which satisfies $P(H) \subset D(S)$, $SP = 0$ as well as $P(D) \subset D(A)$, $AP(D) \subset D(A)$. Moreover, we introduce $P_S : H \to H$ as

$$P_S f = Pf + (f, 1), \quad f \in H.$$

(D6) *The invariant measure.* Let $\mu$ be invariant for $(L, D)$ in the sense that

$$(Lf, 1) = \int_W Lf \, d\mu = 0 \quad \text{for all} \quad f \in D.$$

(D7) *Semigroup conservativity.* $1 \in D(L)$ and $L1 = 0$. Note that this implies $T_t 1 = 1$, as $\frac{d}{dt}T_t 1 = T_t L 1 = 0$ for all $t \in [0, \infty)$.

Now the first three hypocoercivity assumptions read as follows.

**Assumption (A1).** (Algebraic relation) *Assume that* $PAP_{|D} = 0$.

**Assumption (A2).** (Microscopic coercivity) *There exists $\Lambda_m > 0$ such that*

$$-(Sf, f) \geq \Lambda_m \|(I - P_S)f)\|^2 \quad \text{for all} \quad f \in D.$$

**Assumption (A3).** (Macroscopic coercivity) Define $(G, D)$ by $G = PA^2P$ on $D$. Assume that $(G, D)$ is essentially self-adjoint on $H$ (or essentially m-dissipative on $H$ equivalently). Moreover, assume that there exists $\Lambda_M > 0$ such that

$$\|APf\|^2 \geq \Lambda_M \|Pf\|^2 \quad \text{for all} \quad f \in D.$$

In this hypocoercivity setting, one can introduce a suitable bounded linear operator $B$ on $H$ as follows. It is defined as the unique extension of $(B, D((AP)^*))$ to a continuous linear operator on $H$ where

$$B = (I + (AP)^*AP)^{-1}(AP)^* \quad \text{on} \quad D((AP)^*).$$

For the fact that $B$ extends to a bounded linear operator on $H$, consider the original references stated above or see [25, Theo. 5.1.9]. Using the bounded linear operator $B$ one can state the last hypocoercivity assumption.

**Assumption (A4).** (Boundedness of auxiliary operators) The operators $(BS, D)$ and $(BA(I - P), D)$ are bounded and there exist constants $c_1 < \infty$ and $c_2 < \infty$ such that for all $f \in D$

$$\|BSf\| \leq c_1 \|(I - P)f\| \quad \text{and} \quad \|BA(I - P)f\| \leq c_2 \|(I - P)f\|.$$

In applications, in particular for the Langevin dynamics later on, the upcoming lemma is a useful tool to verify (A4).

**Lemma 5.1** *Assume that $(G, D)$ is essentially self-adjoint.*

(i) *If there exists $c_3 \in \mathbb{R}$ such that*

$$\|(BS)^*g\| \leq c_3 \|g\| \quad \text{for all} \quad g = (I - G)f, \quad f \in D,$$

*then the first inequality in (A4) holds with $c_1 = c_3$.*

(ii) *If there exist $c_4 < \infty$ such that*

$$\|(BA)^*g\| \leq c_4 \|g\| \quad \text{for all} \quad g = (I - G)f, \quad f \in D,$$

*then the second inequality in (A4) is satisfied with $c_2 = c_4$.*

Next we formulate the hypocoercivity theorem. The proof can be found in [17], nevertheless we describe in detail how to explicitly compute the constants determining the speed of convergence, compare also [18, Theorem 1.1].

**Theorem 5.2** *Assume that (D) and (A1)–(A4) hold. Then for each $\theta_1 \in (1, \infty)$ there exist $\theta_2 \in (0, \infty)$ such that for each $g \in H$ we have*

$$\|T_t g - (g, 1)\| \leq \theta_1 e^{-\theta_2 t} \|g - (g, 1)\| \quad \text{for all} \quad t \geq 0,$$

*where $(T_t)_{t \geq 0}$ denotes the $C_0$-semigroup introduced in (D2). The constant $\theta_2$ is explicitly computable in terms of $\Lambda_m, \Lambda_M, c_1$ and $c_2$ and given as*

$$\frac{1}{4} \frac{\theta_1 - 1}{\theta_1} \frac{\min\{\Lambda_m, c_1\}}{(1 + c_1 + c_2)\left(1 + \frac{1 + \Lambda_M}{2\Lambda_M}(1 + c_1 + c_2)\right) + \frac{1}{2} \frac{\Lambda_M}{1 + \Lambda_M}} \frac{\Lambda_M}{1 + \Lambda_M}.$$

***Proof*** In view of [17, Theorem 2.18] we first choose $\delta > 0$ such that

$$\frac{\Lambda_M}{1 + \Lambda_M} - (1 + c_1 + c_2)\frac{\delta}{2} > 0.$$

and then $\varepsilon > 0$ small enough such that

$$\Lambda_m - \varepsilon(1 + c_1 + c_2)\left(1 + \frac{1}{2\delta}\right) > 0,$$

as well. This particular choice ensures that

$$\min\left\{\frac{\Lambda_M}{1 + \Lambda_M} - (1 + c_1 + c_2)\frac{\delta}{2}, \Lambda_m - \varepsilon(1 + c_1 + c_2)\left(1 + \frac{1}{2\delta}\right)\right\} > 0.$$

Now chose $\kappa > 0$ smaller or equal than the minimum above. Again, by [17, Theorem 2.18] we obtain

$$\|T_t g - (g, 1)\| \leq \kappa_1 e^{-\kappa_2 t} \|g - (g, 1)\| \quad \text{for all} \quad t \geq 0,$$

for

$$\kappa_1 = \sqrt{\frac{1 + \varepsilon}{1 - \varepsilon}} \quad \text{and} \quad \kappa_2 = \frac{\kappa}{1 + \varepsilon}.$$

To explicitly compute $\theta_1$ and $\theta_2$ as promised in the assertion, we have to specify $\delta > 0$, $\varepsilon > 0$ and the corresponding $\kappa > 0$. We use the strategy from [18, Theorem 1.1]. W.l.o.g. we can assume $\Lambda_m \leq c_1$, since otherwise we can replace $\Lambda_m$ with $\min\{\Lambda_m, c_1\}$ in (A2). We set

$$\delta = \frac{\Lambda_M}{1 + \Lambda_M} \frac{1}{1 + c_1 + c_2}.$$

Moreover, we define

$$r_{\Lambda_M, c_1} = (1 + c_1 + c_2)\left(1 + \frac{1 + \Lambda_M}{2\Lambda_M}(1 + c_1 + c_2)\right) \quad \text{and} \quad s_{\Lambda_M} = \frac{1}{2}\frac{\Lambda_M}{1 + \Lambda_M}.$$

For arbitrary $v \in (0, \infty)$ we chose $\varepsilon = \frac{v}{1+v} \frac{\Lambda_m}{r_{\Lambda_M,c_1}+s_{\Lambda_M}}$. As $\Lambda_m \leq c_1$ one can check that $\varepsilon \in (0, 1)$. Since $\varepsilon(r_{\Lambda_M,c_1} + s_{\Lambda_M}) = \frac{v}{1+v}\Lambda_m < \Lambda_m$ we get

$$\Lambda_m - \varepsilon r_{\Lambda_M,c_1} \geq \varepsilon s_{\Lambda_M} = \frac{v}{1+v} \frac{\Lambda_m}{r_{\Lambda_M,c_1}+s_{\Lambda_M}} s_{\Lambda_M}.$$

In particular $\kappa = \frac{v}{1+v} \frac{\Lambda_m}{r_{\Lambda_M,c_1}+s_{\Lambda_M}} s_{\Lambda_M}$ is a valid choice. The convergence rate in terms of $\kappa_1$ and $\kappa_2$ is given by

$$\kappa_1 = \sqrt{\frac{1+\varepsilon}{1-\varepsilon}} = \sqrt{\frac{1+v+\frac{\Lambda_m}{r_{\Lambda_M,c_1}+s_{\Lambda_M}}v}{1+v-\frac{\Lambda_m}{r_{\Lambda_M,c_1}+s_{\Lambda_M}}v}} \leq \sqrt{1+2v+v^2} = 1+v \quad \text{and}$$

$$\kappa_2 = \frac{\kappa}{1+\varepsilon} > \frac{1}{2}\kappa.$$

Therefore, choosing $\theta_1 = 1+v$ and $\theta_2 = \frac{1}{2}\kappa$ yields the claimed rate of convergence. □

## 6 Hypocoercivity of the langevin dynamics

The aim of this section is to show exponential convergence to equilibrium in our extended hypocoercivity framework. Hypocoercivity refers to the semigroup generated by the infinite-dimensional Langevin operator. At this point we emphasize that the calculations and arguments to verify the data conditions (D1)–(D7) and the hypocoercivity Assumptions (A1)–(A4) below are similar to the associated dual statement in the Fokker plank setting in [13] and including domain issues in [17, 18]. However, the infinite-dimensionality of our problem results in more challenging calculations and the need of advanced arguments. Consequently we take a closer look at the verification of (D1)–(D7) and (A1)–(A4), again.

Throughout this section we assume Assumption 3.4. Hence by Theorem 3.6 the data conditions (D1)–(D4), (D6) and (D7) are already fulfilled. To check the rest we have to specify the orthogonal projection $P_S$.

**Definition 6.1** Define $P_S : L^2(\mu^\Phi) \to L^2(\mu^\Phi)$ by

$$P_S f = \int_V f \, d\mu_2,$$

where the integration is understood w.r.t. the second variable. An application of Fubini's theorem and the fact that $(V, \mathcal{B}(V), \mu_2)$ is a probability space shows that $P_S$ is a well-defined orthogonal projection on $L^2(\mu^\Phi)$ satisfying

$$P_S f \in L^2(\mu_1^\Phi) \quad \text{and} \quad \|P_S f\|_{L^2(\mu_1^\Phi)} = \|P_S f\|_{L^2(\mu^\Phi)}, \quad f \in L^2(\mu^\Phi).$$

In the definition above we canonically embed $L^2(\mu_1^\Phi)$ into $L^2(\mu^\Phi)$. Using that $\mu^\Phi$ is a probability measure one can check that the map $P : L^2(\mu^\Phi) \to L^2(\mu^\Phi)$ given as

$$Pf = P_S f - (f, 1)_{L^2(\mu^\Phi)}, \quad f \in L^2(\mu^\Phi),$$

is an orthogonal projection with

$$Pf \in L^2(\mu_1^\Phi) \quad \text{and} \quad \|Pf\|_{L^2(\mu_1^\Phi)} = \|Pf\|_{L^2(\mu^\Phi)}, \quad f \in L^2(\mu^\Phi).$$

It is important to mention that for each $f \in \mathcal{F}C_b^\infty(B_W)$, the functions $P_S f$ and $Pf$ admit a unique version in $\mathcal{F}C_b^\infty(B_U)$. For notational convenience we write for $f \in \mathcal{F}C_b^\infty(B_W)$

$$f_S = P_S f \quad \text{and} \quad f_P = Pf.$$

The next lemma includes the proof of the data condition (D5) and tells us more about the structure of the infinite-dimensional Langevin operator. Recall the operators defined in Definition 3.9. To avoid to heavy notation we assume w.l.o.g. $\min_{k \in \mathbb{N}}\{m_k \mid n \leq m_k\} = n$ for all $n \in \mathbb{N}$.

**Lemma 6.2** *The data condition (D5) holds true, i.e. $P(L^2(\mu^\Phi)) \subset D(S)$, $SP = 0$ as well as $P(\mathcal{F}C_b^\infty(B_W)) \subset D(A_\Phi)$, $A_\Phi P(\mathcal{F}C_b^\infty(B_W)) \subset D(A_\Phi)$. Moreover, for all $f \in \mathcal{F}C_b^\infty(B_W, n)$, $n \in \mathbb{N}$, we have*

(i) $A_\Phi Pf = -(v, Q_2^{-1} K_{12} D_1 f_S)_V$,
(ii) $A_\Phi^2 Pf = \sum_{i,j=1}^n (v, Q_2^{-1} K_{12} d_i)_V (v, Q_2^{-1} K_{12} d_j)_V \partial_{ij} f_S$
$- (u, Q_1^{-1} C D_1 f_S)_U - (D\Phi(u), C D_1 f_S)_U$ and
(iii) $P_S A_\Phi^2 Pf = \operatorname{tr}[C D_1^2 f_S] - (u, Q_1^{-1} C D_1 f_S)_U - (D\Phi(u), C D_1 f_S)_U.$

**Proof** Given $Pf \in P(L^2(\mu^\Phi))$. Choose a sequence $(f_k)_{k \in \mathbb{N}} \subset \mathcal{F}C_b^\infty(B_W)$ converging to $f$ in $L^2(\mu^\Phi)$. Obviously $(Pf_k)_{k \in \mathbb{N}} \subset \mathcal{F}C_b^\infty(B_W)$ converges to $Pf$ in $L^2(\mu^\Phi)$. For all $k \in \mathbb{N}$, it holds $SPf_k = 0$, since $Pf_k$ only depends on the first variable. Using the fact that $(S, D(S))$ is a closed operator we can conclude that $Pf \in D(S)$ and $SP = 0$.

As $P(\mathcal{F}C_b^\infty(B_W)) \subset \mathcal{F}C_b^\infty(B_W)$, it immediately follows that $P(\mathcal{F}C_b^\infty(B_W)) \subset D(A_\Phi)$. Item (i) follows, since for all $f \in \mathcal{F}C_b^\infty(B_W)$ and $(u, v) \in W$ it holds

$$A_\Phi Pf(u, v) = -(v, Q_2^{-1} K_{12} D_1 f_P(u))_V = -(v, Q_2^{-1} K_{12} D_1 f_S(u))_V.$$

To show $A_\Phi P(\mathcal{F}C_b^\infty(B_W)) \subset D(A_\Phi)$, consider a function $f \in \mathcal{F}C_b^\infty(B_W)$ and a sequence of cut-off functions $(\varphi_m)_{m \in \mathbb{N}} \subset C_0^\infty(\mathbb{R}^n)$ provided by Lemma 3.2. Here, $n \in \mathbb{N}$ corresponds to the $n$ such that $f \in \mathcal{F}C_b^\infty(B_W, n)$. Now we define $\varphi_m^n(v) = \varphi_m(P_n(v))$ and the sequence $(g_m)_{m \in \mathbb{N}}$ by

$$W \ni (u, v) \mapsto g_m(u, v) = \varphi_m^n(v) A_\Phi Pf(u, v) = -\varphi_m^n(v)(v, Q_2^{-1} K_{12} D_1 f_S(u))_V.$$

Hence, the lac of boundedness of $A_\Phi P f$ in the second variable is compensated with the sequence of cut-off functions and therefore $(g_m)_{m\in\mathbb{N}} \subset \mathcal{F}C_b^\infty(B_W)$. Using the product rule we can calculate for all $m \in \mathbb{N}$

$$\begin{aligned}A_\Phi g_m &= (u, Q_1^{-1} K_{21} D_2(A_\Phi P f))_U \varphi_m^n + (D\Phi, K_{21} D_2(A_\Phi P f))_U \varphi_m^n \\ &\quad - (v, Q_2^{-1} K_{12} D_1(A_\Phi P f))_V \varphi_m^n \\ &\quad + \big((D\Phi, K_{21} D_2 \varphi_m^n)_U + (u, Q_1^{-1} K_{21} D_2 \varphi_m^n)_U\big) A_\Phi P f.\end{aligned}$$

The theorem of dominated convergence, implies that $(A_\Phi g_m)_{m\in\mathbb{N}}$ converges to

$$(u, Q_1^{-1} K_{21} D_2(A_\Phi P f))_U + (D\Phi, K_{21} D_2(A_\Phi P f))_U - (v, Q_2^{-1} K_{12} D_1(A_\Phi P f))_V,$$

in $L^2(\mu^\Phi)$. Obviously $(g_m)_{m\in\mathbb{N}}$ converges to $A_\Phi P f$ in $L^2(\mu^\Phi)$, as $m \to \infty$. Since $(A_\Phi, D(A_\Phi))$ is a closed operator we can conclude $A_\Phi P(\mathcal{F}C_b^\infty(B_W)) \subset D(A_\Phi)$ and the function defined right above is $A_\Phi^2 P f$.

In order to show Item $(ii)$ we have to calculate $D_1(A_\Phi P f)$ and $D_2(A_\Phi P f)$. Due to the structure of $A_\Phi P f$, $D_2(A_\Phi P f)$ is easily derived to be $-Q_2^{-1} K_{12} D_1 f_S$. Using the fact that $A_\Phi P f$ only depends on the first $n$-directions in the first variable one can calculate for all $(u, v) \in W$

$$D_1(A_\Phi P f)(u, v) = -\sum_{i,j=1}^n (v, Q_2^{-1} K_{12} d_j)_V \partial_{ij} f_S(u) d_i.$$

To sum up we get

$$\begin{aligned}A_\Phi^2 P f(u, v) &= \sum_{i,j=1}^n (v, Q_2^{-1} K_{12} d_i)_V (v, Q_2^{-1} K_{12} d_j)_V \partial_{ij,1} f_S(u) \\ &\quad - (u, Q_1^{-1} K_{21} Q_2^{-1} K_{12} D_1 f_S(u))_U \\ &\quad - (D\Phi(u), K_{21} Q_2^{-1} K_{12} D_1 f_S(u))_U.\end{aligned}$$

Using Lemma 3.3 we obtain

$$\begin{aligned}\int_V (v, Q_2^{-1} K_{12} d_i)_V (v, Q_2^{-1} K_{12} d_j)_V \, d\mu_2(v) \\ = (Q_2 Q_2^{-1} K_{12} d_i, Q_2^{-1} K_{12} d_j)_V = (K_{12} d_i, Q_2^{-1} K_{12} d_j)_V = (C d_i, d_j)_U,\end{aligned}$$

for all $1 \le i, j \le n$. We can conclude

$$\begin{aligned}&P_S A_\Phi^2 P f(u, v) \\ &= \sum_{i,j=1}^n \partial_{ij} f_S(u) (C d_i, d_j)_U - (u, Q_1^{-1} C D_1 f_S(u))_U - (D\Phi(u) C D_1 f_S(u))_U\end{aligned}$$

$$= \text{tr}[CD_1^2 f_S(u)] - (u, Q_1^{-1} CD_1 f_S(u))_U - (D\Phi(u), CD_1 f_S(u))_U.$$

□

Again, note that the proof above is based on [18, Lemma 3.4] and lifted to our infinite-dimensional setting. Without further tools we can directly show the algebraic relation stated in Sect. 5 Assumption (A1).

**Lemma 6.3** *For all $f \in \mathcal{F}C_b^\infty(B_W)$ it holds $PA_\Phi P f = 0$.*

**Proof** Using Item $(i)$ from Lemma 6.2 it holds for all $f \in \mathcal{F}C_b^\infty(B_W)$ and $(u, v) \in W$

$$A_\Phi P f(u, v) = -(v, Q_2^{-1} K_{12} D_1 f_S(u))_V.$$

Since the Gaussian measure $\mu_2$ is centered we have

$$\int_V (v, x)_V \, d\mu_2(v) = 0 \quad \text{for all} \quad x \in V.$$

Hence we can conclude that $P_S A_\Phi P f = 0$. The fact that

$$(A_\Phi P f, 1)_{L^2(\mu^\Phi)} = (P_S A_\Phi P f, 1)_{L^2(\mu_1^\Phi)} = 0,$$

yields $PA_\Phi P f = 0$ as desired. □

The next checkpoint on our list is Assumption (A2), namely the microscopic coercivity. I.e., we address the symmetric part $S$ of $L_\Phi$ and need an advanced dissipativity statement. In applications this can be checked via the Poincaré inequality for infinite-dimensional Gaussian measures stated in Proposition 3.12. Moreover, the remark below can be useful.

**Remark 6.4** In [11, Section 10.5.2.] it is shown that Item $(iii)$ of Assumption 3.10 holds true, if there exist constants $M_1, \tilde{\omega}_1 \in (0, \infty)$ such that

$$\|K_{22}^{\frac{1}{2}} \exp(-t K_{22} Q_2^{-1}) K_{22}^{-\frac{1}{2}}\|_{\mathcal{L}(V)} \le M_1 \exp(-\tilde{\omega}_1 t), \quad t \in [0, \infty).$$

In this case $\omega_1 = \frac{\tilde{\omega}_1}{M_1}$.

**Lemma 6.5** *Assume Item $(iii)$ from Assumption 3.10 is valid. Then it holds*

$$-(Sf, f)_{L^2(\mu^\Phi)} \ge \omega_1 \|(I - P_S)f\|^2_{L^2(\mu^\Phi)} \quad \text{for all} \quad f \in \mathcal{F}C_b^\infty(B_W),$$

*i.e. (A2) with $\Lambda_m = \omega_1$.*

**Proof** In view of the Poincaré inequality from Assumption 3.10 and the integration by parts formula from Theorem 3.6, we obtain for all $f \in \mathcal{F}C_b^\infty(B_W)$

$$-(Sf, f)_{L^2(\mu^\Phi)} = \int_U \int_V (K_{22} D_2 f, D_2 f)_V \, d\mu_2 d\mu_1^\Phi$$

$$\geq \omega_1 \int_U \|f(u,\cdot) - (f(u,\cdot), 1)_{L^2(\mu_2)}\|^2_{L^2(\mu_2)} \, d\mu_1^\Phi(u)$$
$$= \omega_1 \|(I - P_S)f\|^2_{L^2(\mu^\Phi)}.$$

□

The upcoming statements are devoted to prove Assumption (A3), i.e. the macroscopic coercivity. Recall the operator $G$ from Sect. 5.

**Proposition 6.6** *Assume Item (i) from Assumption 3.10 is valid. Then the linear operator* $(G, \mathcal{F}C_b^\infty(B_W)) = (PA_\Phi^2 P, \mathcal{F}C_b^\infty(B_W))$ *is essentially self-adjoint in* $L^2(\mu^\Phi)$. *Moreover, for* $f \in \mathcal{F}C_b^\infty(B_W)$ *and* $(u, v) \in W$ *it is given via the formula*

$$Gf(u, v) = \mathrm{tr}[CD_1^2 f_S(u)] - (u, Q_1^{-1} CD_1 f_S(u))_U - (D\Phi(u), CD_1 f_S(u))_U.$$

*Proof* For the moment we consider the operator $(N, \mathcal{F}C_b^\infty(B_U))$, which is for $f \in \mathcal{F}C_b^\infty(B_U)$ and $u \in U$ defined by

$$Nf(u) = \mathrm{tr}[CD^2 f(u)] - (u, Q_1^{-1} CDf(u))_U - (D\Phi(u), CDf(u))_U.$$

The integration by parts formula from Theorem 3.11 implies for all $f \in \mathcal{F}C_b^\infty(B_W)$

$$(A_\Phi^2 Pf, 1)_{L^2(\mu^\Phi)} = (P_S A_\Phi^2 Pf, 1)_{L^2(\mu_1^\Phi)} = (Nf_S, 1)_{L^2(\mu_1^\Phi)} = 0.$$

This implies

$$Gf = PA_\Phi^2 Pf = P_S A_\Phi^2 Pf - (A_\Phi^2 Pf, 1)_{L^2(\mu^\Phi)} = P_S A_\Phi^2 Pf.$$

By Lemma 6.2 Item $(iii)$ we obtain the representation of $(G, \mathcal{F}C_b^\infty(B_W))$ promised in the assertion. Note that the calculation above also shows $Gf = Nf_S$. Another application of the integration by parts formula from Theorem 3.11 shows that $(G, \mathcal{F}C_b^\infty(B_W))$ is symmetric and dissipative in $L^2(\mu^\Phi)$. As densely defined symmetric operators on a Hilbert space are essentially self-adjoint if and only if they are essentially m-dissipative, it is left to show that $(Id - G)(\mathcal{F}C_b^\infty(B_W))$ is dense in $L^2(\mu^\Phi)$. We proof this by showing that

$$((I - G)f, g)_{L^2(\mu^\Phi)} = 0 \quad \text{for all} \quad f \in \mathcal{F}C_b^\infty(B_W), \tag{7}$$

implies $g = 0$. So suppose the statement (7) is true. Given $n \in \mathbb{N}$ and an arbitrary element $f \in \mathcal{F}C_b^\infty(B_U, n)$. Choose a sequence of cut-off function $(\varphi_m)_{m \in \mathbb{N}} \subset C_0^\infty(\mathbb{R}^n)$ provided by Lemma 3.2. The sequence $(f_m)_{m \in \mathbb{N}}$ defined by $W \ni (u, v) \mapsto f_m(u, v) = f(u)\varphi_m(P_n(v)) \in \mathbb{R}$ is in $\mathcal{F}C_b^\infty(B_W)$ and for all $m \in \mathbb{N}$ it holds

$$0 = ((I - G)f_m, g)_{L^2(\mu^\Phi)}$$
$$= (f_m, g)_{L^2(\mu^\Phi)} - \int_W \varphi_m(P_n(v)) Nf(u) g(u, v) \, d\mu^\Phi(u, v)$$

$$\to (f, g_S)_{L^2(\mu_1^\Phi)} - (Nf, g_S)_{L^2(\mu_1^\Phi)},$$

as $m \to \infty$ by the theorem of dominated convergence. Therefore

$$0 = ((I - N)f, g_S)_{L^2(\mu_1^\Phi)} = 0 \quad \text{for all} \quad f \in \mathcal{F}C_b^\infty(B_U, n).$$

Since $n \in \mathbb{N}$ was arbitrary it holds

$$0 = ((I - N)f, g_S)_{L^2(\mu_1^\Phi)} = 0 \quad \text{for all} \quad f \in \mathcal{F}C_b^\infty(B_U).$$

Since $(Id - N)(\mathcal{F}C_b^\infty(B_U)$ is dense in $L^2(\mu_1^\Phi)$ by Theorem 3.11, it follows $g_S = 0$ in $L^2(\mu_1^\Phi)$. Thus for all $f \in \mathcal{F}C_b^\infty(B_W)$

$$(f, g)_{L^2(\mu^\Phi)} = (Gf, g)_{L^2(\mu^\Phi)} = (Nf_S, g_S)_{L^2(\mu_1^\Phi)} = 0.$$

Therefore $g = 0$ by the density of $\mathcal{F}C_b^\infty(B_W)$ in $L^2(\mu^\Phi)$. □

**Lemma 6.7** *Assume that Item (iv) from Assumption 3.10 holds true. Then for all $f \in \mathcal{F}C_b^\infty(B_W)$ we have*

$$\|A_\Phi Pf\|^2_{L^2(\mu^\Phi)} \geq \omega_2 \|Pf\|^2_{L^2(\mu^\Phi)},$$

*i.e. Assumption (A3) with $\Lambda_M = \omega_2$.*

*Proof* Given $f \in \mathcal{F}C_b^\infty(B_W)$. It holds by Lemma 3.3 and Item (iv) from Assumption 3.10

$$\begin{aligned}\|A_\Phi Pf\|^2_{L^2(\mu^\Phi)} &= \int_U \int_V (v, Q_2^{-1} K_{12} D_1 f_S(u))_V^2 \, d\mu_2(v) d\mu_1^\Phi(u) \\ &= \int_U (Q_2^{-1} K_{12} Df_S(u), K_{12} Df_S(u))_V \, d\mu_1^\Phi(u) \\ &\geq \omega_2 \int_U \left(f_S(u) - (f_S, 1)_{L^2(\mu_1^\Phi)}\right)^2 d\mu_1^\Phi(u) \\ &= \omega_2 \|Pf\|^2_{L^2(\mu^\Phi)}.\end{aligned}$$

□

With the help of Lemma 5.1 we want verify Assumption (A4), which addresses the boundedness of the auxiliary operators. The central tool to apply the lemma is the elliptic a priori estimate of Dolbeaut, Mouhot and Schmeiser from [13, Sec. 2, Eq. (2.2), Lem 8]. Since Dolbeaut, Mouhot and Schmeiser only deals with finite dimensions in [13], we had to lift their result to our infinite-dimensional setting. Recall that this lifting is discussed in Sect. 3.

**Proposition 6.8** *Assume that Item (v) from Assumption* 3.10 *holds true with constant* $C_1 \in (0, \infty)$. *Then the operator* $(BS, \mathcal{F}C_b^\infty(B_W))$ *is bounded with*

$$\|BSf\|_{L^2(\mu^\Phi)} \leq C_1 \|(I - P)f\|_{L^2(\mu^\Phi)} \quad \text{for all} \quad f \in \mathcal{F}C_b^\infty(B_W),$$

*i.e. the first inequality in (A4) holds true with* $c_1 = C_1$.

**Proof** To proof the statement we aim to use Lemma 5.1. Given $f \in \mathcal{F}C_b^\infty(B_W)$ and set $g = (Id - G)f$. As shown in [17, Proposition 2.15] it holds

$$B^* g(u, v) = A_\Phi P f(u, v) = -(v, Q_2^{-1} K_{12} D_1 f_S(u))_V,$$

where the last equality is due to Lemma 6.3. An approximation with cut-off functions, as in the proof of Lemma 6.2, shows that $A_\Phi P f \in D(S)$ with

$$(BS)^* g(u, v) = S^* A_\Phi P f(u, v) = S A_\Phi P f(u, v) = (v, Q_2^{-1} K_{22} Q_2^{-1} K_{12} D_1 f_S(u))_V.$$

Using Lemma 3.3 and Assumption 3.10 we get

$$\begin{aligned}\|(BS)^* g\|^2_{L^2(\mu^\Phi)} &= \int_U \int_V (v, Q_2^{-1} K_{22} Q_2^{-1} K_{12} D_1 f_S)^2_V \, d\mu_2 d\mu_1^\Phi \\ &= \int_U (K_{22} Q_2^{-1} K_{12} D_1 f_S(u), Q_2^{-1} K_{22} Q_2^{-1} K_{12} D_1 f_S)_V \, d\mu_1^\Phi \\ &\leq C_1^2 \|g\|^2_{L^2(\mu^\Phi)}.\end{aligned}$$

Therefore the claim follows by Lemma 5.1. □

**Proposition 6.9** *Assume that Item (i) and (ii) from Assumption* 3.10 *is satisfied. Then the operator* $(BA_\Phi(I - P), \mathcal{F}C_b^\infty(B_W))$ *is bounded and*

$$\|BA_\Phi(I - P)f\|_{L^2(\mu^\Phi)} \leq 2\sqrt{2} \|(I - P)f\|_{L^2(\mu^\Phi)} \quad \text{for all} \quad f \in \mathcal{F}C_b^\infty(B_W).$$

*I.e., the second inequality in (A4) holds true with* $c_2 = 2\sqrt{2}$.

**Proof** As in the proposition above we want to apply Lemma 5.1. For $f \in \mathcal{F}C_b^\infty(B_W)$ and $g = (I - G)f$ it holds by [17, Proposition 2.15]

$$(BA_\Phi)^* f = -A_\Phi^2 P f.$$

Using Lemma 6.3 we can calculate

$$\|A_\Phi^2 P f\|^2_{L^2(\mu^\Phi)} = \int_U \int_V \left( \sum_{i,j=1}^n (v, Q_2^{-1} K_{12} d_i)_V (v, Q_2^{-1} K_{12} d_j)_V \partial_{ij} f_S(u) \right)^2 d\mu_2(v) d\mu_1^\Phi(u)$$

$$-2\sum_{i,j=1}^{n}\int_{V}(v,Q_2^{-1}K_{12}d_i)_V(v,Q_2^{-1}K_{12}d_j)_V\,\mathrm{d}\mu_2(v)$$
$$\times\int_{U}\partial_{ij}f_S(u)\Big((u,Q_1^{-1}CD_1f_S(u))_U+(D\Phi(u),CD_1f_S(u))_U\Big)\mathrm{d}\mu_1^{\Phi}(u)$$
$$+\int_{U}\Big((u,Q_1^{-1}CD_1f_S(u))_U+(D\Phi(u),CD_1f_S(u))_U\Big)^2\mathrm{d}\mu_1^{\Phi}(u).$$

Moreover, by Lemma 3.3 it holds for all $i,j,k,l\in\mathbb{N}$

$$\int_V (v,Q_2^{-1}K_{12}d_i)_V(v,Q_2^{-1}K_{12}d_j)_V\,\mathrm{d}\mu_2(v)=c_{ij}$$
$$\int_V (v,Q_2^{-1}K_{12}d_i)_V(v,Q_2^{-1}K_{12}d_j)_V(v,Q_2^{-1}K_{12}d_k)_V(v,Q_2^{-1}K_{12}d_l)_V\,\mathrm{d}\mu_2(v)$$
$$=c_{ij}c_{kl}+c_{ik}c_{jl}+c_{il}c_{jk},$$

where $c_{ij}=(Cd_i,d_j)_U$. Recall the definition of $N$ from Proposition 6.6. Substituting the second equation into the first, one arrives at

$$\|A_\Phi^2 Pf\|^2_{L^2(\mu^\Phi)}=\int_U \mathrm{tr}[CD_1^2 f_S(u)]^2-2\mathrm{tr}[CD_1^2 f_S(u)]((u,Q_1^{-1}CD_1f_S(u))_U$$
$$+(D\Phi(u),CD_1f_S(u))_U)+((u,Q_1^{-1}CD_1f_S(u))_U$$
$$+(D\Phi(u),CD_1f_S(u))_U)^2+2\mathrm{tr}[(CD^2 f_S(u))^2]\mathrm{d}\mu_1^\Phi(u)$$
$$=\int_U Nf_S(u)^2\mathrm{d}\mu_1^\Phi(u)+2\int_U \mathrm{tr}[(CD_1^2 f_S(u))^2]\mathrm{d}\mu_1^\Phi(u).$$

Since $Pg=Pf-PGf=Pf-NPf$, we get by the regularity estimates from Theorem 3.11

$$\|A_\Phi^2 Pf\|^2_{L^2(\mu^\Phi)}=\int_U NPf(u)^2\mathrm{d}\mu_1^\Phi(u)+2\int_U \mathrm{tr}[(CD^2 Pf(u))^2]\mathrm{d}\mu_1^\Phi(u)$$
$$\leq 2\int_U Pf(u)^2+Pg(u)^2\mathrm{d}\mu_1^\Phi(u)+4\int_U Pg(u)^2\mathrm{d}\mu_1^\Phi(u)$$
$$\leq 8\|Pg\|^2_{L^2(\mu^\Phi)}\leq 8\|g\|^2_{L^2(\mu^\Phi)}.$$

Hence the claim follows by Lemma 5.1 $(ii)$. □

In conclusion, we can apply Theorem 5.2 to get the promised hypocoercivity result, summarized in the next theorem.

**Theorem 6.10** *Assume that Assumption 3.4 and 3.10 hold true. Then the operator $(L_\Phi,\mathcal{F}C_b^\infty(B_W))$ is essentially m-dissipative in $L^2(\mu^\Phi)$. The semigroup $(T_t)_{t\geq 0}$ generated by the closure of $(L_\Phi,\mathcal{F}C_b^\infty(B_W))$ is hypocoercive. In particular, for each $\theta_1\in(1,\infty)$ there is some $\theta_2\in(0,\infty)$ such that*

$$\|T_tg-(g,1)_{L^2(\mu^\Phi)}\|_{L^2(\mu^\Phi)}\leq\theta_1 e^{-\theta_2 t}\|g-(g,1)_{L^2(\mu^\Phi)}\|_{L^2(\mu^\Phi)}\quad\text{for all}\ \ t\geq 0.$$

*The constant $\theta_2$ determining the speed of convergence can be explicitly computed as*

$$\theta_2 = \frac{1}{2}\frac{\theta_1-1}{\theta_1}\frac{\min\{\omega_1,C_1\}}{(1+2\sqrt{2}+C_1)\left(1+\frac{1+\omega_2}{2\omega_2}(1+2\sqrt{2}+C_1)\right)+\frac{1}{2}\frac{\omega_2}{1+\omega_2}}\frac{\omega_2}{1+\omega_2}.$$

*Proof* By the considerations above we can determine the constants from (A2)-(A4). I.e., it holds

$$\Lambda_M = \omega_2 \quad \Lambda_m = \omega_1 \quad c_1 = C_1 \quad \text{and} \quad c_2 = 2\sqrt{2}.$$

Applying Theorem 5.2 with these constants yields the claim. □

**Remark 6.11** In view of Remark 3.12 from [18] one can study the rate of convergence in terms of $K_{22}$. Indeed, if we assume that $K_{22} = \alpha Q_2$ for some $\alpha \in (0,\infty)$ one can calculate that $\omega_1 = \alpha\lambda_1$ and $C_1 = \alpha$ are valid constants. Hence for small $\alpha$ we obtain a bad exponential convergence rate. This reflects the observation that for small $\alpha$ the equation becomes almost deterministic.

Since $\theta_2$ converges to zero for $\alpha \to \infty$ the rate of exponential convergence gets arbitrary bad by considering large $\alpha$. This is consistent with the calculations in [18, Remark 3.12], where a finite-dimensional Langevin equation is considered. Indeed in [18, Remark 3.12] the authors show that the exponential convergence rate in a large damping regime, which corresponds to large $\alpha$, is of order $\frac{1}{\alpha}$ and hence arbitrary bad.

We want to end this section with an $L^2$-exponential ergodicity result for the weak solution provided by Theorem 4.7, which is shown in a manifold setting in [24, Corollary 5.2].

**Corollary 6.12** *Assume that Assumption 3.4, 3.8 and 3.10 hold true. Let $\theta_1 \in (1,\infty)$ and $\theta_2 \in (0,\infty)$ be the constants determined by Theorem 6.10 and*

$$X = (\Omega, \mathcal{F}, (\mathcal{F}_t)_{t\geq 0}, (X_t)_{t\geq 0}, (\theta_t)_{t\geq 0}, (\mathbb{P}^w)_{w\in W}),$$

*the $\mu^\Phi$-standard right process providing a weak solution as described in Theorem 4.7. Then for all $t \in (0,\infty)$ and $g \in L^2(\mu^\Phi)$, $T_t g$ is a $\mu^\Phi$-version of $p_t g$. Furthermore, it holds*

$$\left\|\frac{1}{t}\int_{[0,t)} g(X_s)d\lambda(s) - (g,1)_{L^2(\mu^\Phi)}\right\|_{L^2(\mathbb{P}_{\mu^\Phi})}$$
$$\leq \frac{1}{\sqrt{t}}\sqrt{\frac{2\theta_1}{\theta_2}\left(1-\frac{1}{t\theta_2}(1-e^{-t\theta_2})\right)}\|g-(g,1)_{L^2(\mu^\Phi)}\|_{L^2(\mu^\Phi)},$$

*for all $t \in (0,\infty)$. We call a weak solutions $X$ with this property $L^2$-exponentially ergodic, i.e. ergodic with a rate that corresponds to exponential convergence of the corresponding semigroup.*

**Proof** The relation between $T_t g$ and $p_t g$, for $t > 0$ and $g \in L^2(\mu^\Phi)$, is part of Theorem 4.3. To show ergodicity, let $t \in (0, \infty)$ and $g \in L^2(\mu^\Phi)$ be given. For $f = g - (g, 1)_{L^2(\mu^\Phi)}$ it holds

$$\left\| \frac{1}{t} \int_{[0,t)} f(X_s) d\lambda(s) \right\|^2_{L^2(\mathbb{P}_{\mu^\Phi})} = \frac{2}{t^2} \int_{[0,t)} \int_{[0,s)} (T_{s-u} f, f)_{L^2(\mu^\Phi)} d\lambda(u) d\lambda(s)$$

$$\leq \frac{2\|f\|^2_{L^2(\mu^\Phi)}}{t^2} \int_{[0,t)} \int_{[0,s)} \theta_1 e^{-(s-u)\theta_2} d\lambda(u) d\lambda(s)$$

$$= \frac{1}{t} \frac{2\theta_1}{\theta_2} \left(1 - \frac{1}{t\theta_2}(1 - e^{-t\theta_2})\right) \|f\|^2_{L^2(\mu^\Phi)}.$$

Note, to obtain the first equality we argued as in [24, Corollary 5.2]. Afterwards, we used the Cauchy-Schwarz inequality and the hypocoercivity of the semigroup. In the last line we computed the integral. □

**Remark 6.13** From Corollary 6.12 above we can follow, roughly speaking, that time average converges to space average in $L^2(\mathbb{P}_{\mu^\Phi})$ with rate $t^{-\frac{1}{2}}$. If the spectrum of $(L_\Phi, D(L_\Phi))$ contains a negative eigenvalue $-\kappa$ with corresponding eigenvector $g$, then this rate is optimal. Indeed, by a similar reasoning as in the calculation above we get for all $t \in (0, \infty)$

$$\left\| \frac{1}{t} \int_{[0,t)} g(X_s) d\lambda(s) \right\|_{L^2(\mathbb{P}_{\mu^\Phi})} = \frac{1}{\sqrt{t}} \sqrt{\frac{2}{\kappa}\left(1 - \frac{1}{t\kappa}(1 - e^{-t\kappa})\right)} \|g\|_{L^2(\mu^\Phi)}.$$

Equality above holds, as the application of the Cauchy-Schwarz inequality is not necessary. Moreover, note that $(g, 1)_{L^2(\mu^\Phi)} = \frac{1}{-\kappa}(L_\Phi g, 1)_{L^2(\mu^\Phi)} = 0$.

## 7 Examples

We fix a real separable Hilbert space $U$ and a self-adjoint operator $Q \in \mathcal{L}^+(U)$ with corresponding basis of eigenvectors $(e_i)_{i \in \mathbb{N}}$ and a decreasing sequence eigenvalues of $(\lambda_i)_{i \in \mathbb{N}} \subset (0, \infty)$. Moreover, we set $W = U \times U$. For $\alpha_1, \alpha_2 \in (0, \infty)$ such that

$$(\lambda_i^{\alpha_1})_{i \in \mathbb{N}}, (\lambda_i^{\alpha_2})_{i \in \mathbb{N}} \in l^1(\mathbb{N}), \tag{8}$$

we consider $Q_1 = Q^{\alpha_1}$ and $Q_2 = Q^{\alpha_2}$ as covariance operators of two centered Gaussian measures $\mu_1$ and $\mu_2$, respectively. Now we fix a potential $\Phi \in W^{1,2}(U, \mu_1, \mathbb{R})$, which is bounded from below, with bounded gradient and such that $\int_U e^{-\Phi} d\mu_1 = 1$. For $\beta_1, \beta_2 \in [0, \infty)$ we set $K_{12} = Q^{\beta_1}$ and $K_{22} = Q^{\beta_2}$. Since $K_{21} = K_{12}^*$, also $K_{21} = Q^{\beta_1}$. For this particular choice of operators the infinite-dimensional degener-ate stochastic differential equation (1) reads as follows

$$dU_t = Q^{\beta_1-\alpha_2} V_t dt \qquad (9)$$
$$dV_t = -(Q^{\beta_2-\alpha_2} V_t + Q^{\beta_1-\alpha_1} U_t + Q^{\beta_1} D\Phi(U_t))dt + \sqrt{2Q^{\beta_2}} dW_t.$$

The associated generator is given by

$$L_\Phi f = \mathrm{tr}[Q^{\beta_2} D_2^2 f] - (v, Q^{\beta_2-\alpha_2} D_2 f)_U,$$
$$- (u, Q^{\beta_1-\alpha_1} D_2 f)_U - (D\Phi(u), Q^{\beta_1} D_2 f)_U + (v, Q^{\beta_1-\alpha_2} D_1 f)_U.$$

Assuming

$$\beta_2 \leq 2\beta_1, \qquad (10)$$

Item $(iv)$ from Assumption 3.4 is easily verified. Since by construction the other items from Assumption 3.4 are valid, Theorem 3.6 is applicable. I.e., we obtain essential m-dissipativity of $(L_\Phi, \mathcal{F}C_b^\infty(B_W))$ in $L^2(W, \mu^\Phi, \mathbb{R})$.

To construct the corresponding diffusion process $X$ with infinite lifetime, we have to verify the items from Assumption 3.8. In terms of $\alpha_1, \alpha_2, \beta_1$ and $\beta_2$ we need:

(i)
$$(\lambda_i^{\beta_2})_{i\in\mathbb{N}} \in l^1(\mathbb{N}), \qquad (11)$$

to ensure that $K_{22}$ is of trace class.

(ii) A function $\rho \in L^1(\mu^\Phi)$, such that for all $n \in \mathbb{N}$, the function $\rho_n$ defined by

$$W \ni (u, v) \to \rho_n(u, v) = \rho(\overline{P}_n(u), \overline{P}_n(v)) \in \mathbb{R},$$

is in $L^1(\mu^\Phi)$. Furthermore, for all $(u, v) \in W$ we have to check

$$(\overline{P}_n(u), Q^{-\alpha_1+\beta_1}\overline{P}_n(v))_U - (\overline{P}_n(v), Q^{-\alpha_2+\beta_1}\overline{P}_n(u))_U \leq \rho_n(u, v),$$

and $\lim_{n\to\infty} \rho_n = \rho$ in $L^1(\mu^\Phi)$.

If
$$\alpha_1 = \alpha_2, \qquad (12)$$

$\rho = 0$ is the obvious choice. If $\alpha_1$ is not equal to $\alpha_2$, chose $\alpha_1, \alpha_2, \beta_1$ such that $(\overline{P}_n(u), Q^{-\alpha_1+\beta_1}\overline{P}_n(v))_U$ and $(\overline{P}_n(v), Q^{-\alpha_2+\beta_1}\overline{P}_n(u))_U$ define Cauchy sequences in $L^1(\mu^\Phi)$. I.e., assume that

$$(\lambda_i^{\beta_1+\frac{\alpha_2}{2}-\frac{\alpha_1}{2}})_{i\in\mathbb{N}}, (\lambda_i^{\beta_1+\frac{\alpha_1}{2}-\frac{\alpha_2}{2}})_{i\in\mathbb{N}} \in l^1(\mathbb{N}). \qquad (13)$$

This is sufficient, since we obtain for $n \geq m$ and $\gamma \in \{\alpha_1, \alpha_2\}$

$$\int_W |(\overline{P}_n(u), Q_1^{-\gamma+\beta_1}\overline{P}_n(v))_U - (\overline{P}_m(u), Q_1^{-\gamma+\beta_1}\overline{P}_m(v))_U| d\mu^\Phi(u,v)$$

$$\leq \sum_{i=m+1}^n \int_W |\lambda_i^{-\gamma+\beta_1}(u,e_i)_U(v,e_i)_U|e^{-\Phi(u)} d\mu(u,v)$$

$$\leq \sum_{i=m+1}^n \lambda_i^{-\gamma+\beta_1} \int_U |(u,e_i)_U| d\mu_1(u) \int_U |(v,e_i)_U| d\mu_2(v)$$

$$= \sum_{i=m+1}^n \lambda_i^{-\gamma+\beta_1} \int_{\mathbb{R}} |x| d\mathcal{N}(0, \lambda_i^{\alpha_1})(x) \int_{\mathbb{R}} |x| d\mathcal{N}(0, \lambda_i^{\alpha_2})(x)$$

$$= \frac{2}{\pi} \sum_{i=m+1}^n \lambda_i^{-\gamma+\beta_1+\frac{\alpha_1}{2}+\frac{\alpha_2}{2}},$$

where $\mathcal{N}(0, \lambda_i^{\alpha_2})$ and $\mathcal{N}(0, \lambda_i^{\alpha_1})$ are one-dimensional centered Gaussian measures with variances $\lambda_i^{\alpha_2}$ and $\lambda_i^{\alpha_1}$, respectively. Hence, defining $\rho$ as the $L^1(W, \mu^\Phi, \mathbb{R})$ limit of

$$(\overline{P}_n(u), Q^{-\alpha_1+\beta_1}\overline{P}_n(v))_U - (\overline{P}_n(v), Q^{-\alpha_2+\beta_1}\overline{P}_n(u))_U,$$

is another possibility to verify Item (ii) from Assumption 3.8.

To show hypocoercivtiy, with explicitly computable constants, of the semigroup $(T_t)_{t\geq 0}$ generated by $(L_\Phi, D(L_\Phi))$ we use the results from Sect. 6. It is easy to check that the data conditions and the algebraic relation $PAP = 0$ on $D$ are fulfilled. Hence we are left to verify the microscopic/macroscopic coercivity and the boundedness of the auxiliary operators. The operators $(C, D(C))$ and $(Q_1^{-1}C, D(Q_1^{-1}C))$ on $U$ are given by

$$C = K_{21}Q_2^{-1}K_{12} = Q^{-\alpha_2+2\beta_1} \quad \text{and} \quad Q_1^{-1}C = Q^{-\alpha_1-\alpha_2+2\beta_1},$$

respectively. In order to obtain the microscopic coercivity, we use Remark 6.4. I.e., it is enough to assume that

$$\beta_2 - \alpha_2 \leq 0.$$

To get the macroscopic hypocoercivity we need

$$\|D\Phi\|_{L^\infty(\mu_1)} \leq \frac{1}{2\sqrt{\lambda_1}}. \tag{14}$$

Moreover, we have to find a constant $\omega_2 \in (0, \infty)$ s.t. for all $f \in \mathcal{F}C_b^\infty(B_U)$

$$\int_U (Q_2^{-1}K_{12}Df, K_{12}Df)_U d\mu_1^\Phi \geq \omega_2 \int_U (f - (f,1)_{L^2(\mu^\Phi)})^2 d\mu_1^\Phi.$$

If $\Phi : U \to (-\infty, \infty]$ is convex and lower semicontinuous we get in view of Proposition 3.12 for all $f \in \mathcal{F}C_b^\infty(B_U)$

$$\int_U (Q_1 Df, Df)_U \, d\mu_1^\Phi \geq \lambda_1 \int_U \left(f - (f, 1)_{L^2(\mu^\Phi)}\right)^2 d\mu_1^\Phi.$$

In particular, in terms of $\alpha_1, \alpha_2, \beta_1$ we need

$$2\beta_1 - \alpha_2 \leq \alpha_1,$$

to get the macroscopic coercivity with $\omega_2 = \lambda_1$. To obtain boundedness of the auxiliary operators, we have to find a constant $C_1 \in (0, \infty)$ such that for all $f \in \mathcal{F}C_b^\infty(B_W)$ and $g = (I - PA_\Phi^2 P)f$ it holds

$$\int_U (K_{22} Q_2^{-1} K_{12} D_1 fS, Q_2^{-1} K_{22} Q_2^{-1} K_{12} D_1 fS)_U \, d\mu_1^\Phi \leq C_1 \|g\|_{L^2(\mu^\Phi)}^2.$$

By the regularity estimate from Theorem 3.11 we have

$$\int_U (CD_1 fS, Q_1^{-1} CD_1 fS)_U \, d\mu_1^\Phi \leq 4\|g\|_{L^2(\mu^\Phi)}^2.$$

I.e., in terms of $\alpha_1, \alpha_2, \beta_1, \beta_2$ this translates to

$$2\beta_1 - \alpha_1 \leq 2\beta_2 - \alpha_2.$$

In summary, to obtain hypocoercivity we need the following three conditions

$$\beta_2 - \alpha_2 \leq 0, \quad 2\beta_1 - \alpha_1 \leq \alpha_2 \quad \text{and} \quad 2\beta_1 - \alpha_1 \leq 2\beta_2 - \alpha_2.$$

Since the second inequality is implied by the first and the third, we only need

$$\beta_2 - \alpha_2 \leq 0 \quad \text{and} \quad 2\beta_1 - \alpha_1 \leq 2\beta_2 - \alpha_2. \tag{15}$$

In the corollary below we summarize the results we just derived.

**Corollary 7.1** *Assume we are in the setting as described above, with parameters $\alpha_1, \alpha_2 \in (0, \infty)$ such that (8) holds. Moreover, let $\beta_1, \beta_2 \in [0, \infty)$ such that (10) holds true. Then the operator $(L_\Phi, \mathcal{F}C_b^\infty(B_W))$ is essentially m-dissipative in $L^2(W, \mu^\Phi, \mathbb{R})$.*

(i) *If (11) and either (12) or (13) hold true, then there exist a corresponding diffusion process X with infinite lifetime providing a weak solution to (9).*
(ii) *If $\Phi$ is convex and lower-semicontinuous, (14) and (15) are valid, then the strongly continuous contraction semigroup $(T_t)_{t \geq 0}$ corresponding to $(L_\Phi, \mathcal{F}C_b^\infty(B_W))$ is hypocoercive, with explicitly computable constants.*

*If the assumptions from Item (i) and (ii) are valid, then X is an $L^2$-exponentially ergodic weak solution.*

For the rest of this section we want to analyze degenerate stochastic reaction–diffusion equations in the context of the previous example. Non-degenerate equations of this type have already been studied in [9, Section 5].

We set $U = L^2((0, 1), \lambda, \mathbb{R})$ and denote by $(-\Delta, D(\Delta))$ the negative Dirichlet Laplacian, i.e.

$$D(\Delta) = W_0^{1,2}(0, 1) \cap W^{2,2}(0, 1) \subset L^2((0, 1), \lambda, \mathbb{R}),$$
$$-\Delta x = -x''.$$

Moreover, we consider the linear continuous positive self-adjoint operator

$$Q = -\Delta^{-1} : L^2((0, 1), \lambda, \mathbb{R}) \to D(\Delta).$$

As above we define $Q_1, Q_2, K_{22}, K_{21}$ and $K_{12}$ as powers of $Q$ in terms of $\alpha_1, \alpha_2, \beta_1$ and $\beta_2$. Note that $(\sqrt{2}\sin(k\pi \cdot))_{k \in \mathbb{N}}$ is the orthonormal basis of $L^2((0, 1), \lambda, \mathbb{R})$ diagonalizing $Q$ with corresponding eigenvalues $(\lambda_k)_{k \in \mathbb{N}} = (\frac{1}{k^2\pi^2})_{k \in \mathbb{N}}$. Since $(\lambda_k)_{k \in \mathbb{N}} \in l^\theta(\mathbb{N})$ for $\theta \in (\frac{1}{2}, \infty)$ we need $\alpha_1, \alpha_2 \in (\frac{1}{2}, \infty)$, by (8). Additionally we fix a continuous differentiable convex function $\phi : \mathbb{R} \to \mathbb{R}$, which is bounded from below. Assume that there is a constant $C \in [0, \infty)$ such that

$$\sup_{t \in \mathbb{R}} |\phi'(t)| \leq C.$$

For such $\phi$ we consider potentials $\Phi : L^2((0, 1), \lambda, \mathbb{R}) \to \mathbb{R}$ defined by

$$\Phi(u) = \int_{(0,1)} \phi(u(\xi)) d\lambda(\xi).$$

Note the boundedness of $\phi'$ implies that $\phi$ grows at most linear. Moreover, such potentials are lower semicontinuous by Fatou's lemma and bounded from below. Using [9, Proposition 5.2] we know that $\Phi$ is in $W^{1,2}(U, \mu_1, \mathbb{R})$ with $D\Phi(u) = \phi' \circ u$ for $u \in L^2((0, 1), \lambda, \mathbb{R})$. In particular, we obtain

$$\|D\Phi\|_{L^\infty(\mu_1)} = \sup_{t \in \mathbb{R}} |\phi'(t)| \leq C.$$

As $\phi$ is convex the same holds true for $\Phi$. Now the degenerate stochastic reaction–diffusion equation reads as

$$dU_t = (-\Delta)^{-\beta_1 + \alpha_2} V_t dt$$
$$dV_t = -\left((-\Delta)^{-\beta_2 + \alpha_2} V_t + (-\Delta)^{-\beta_1 + \alpha_1} U_t + (-\Delta)^{-\beta_1} \phi' \circ U_t\right) dt$$
$$+ \sqrt{2(-\Delta)^{-\beta_2}} dW_t.$$

The corresponding generator is given by

$$L_\Phi f = \text{tr}[(-\Delta)^{-\beta_2} D_2^2 f] + (v, (-\Delta)^{-\beta_2+\alpha_2} D_2 f)_{L^2(\lambda)}$$
$$+ (u, (-\Delta)^{-\beta_1+\alpha_2} D_2 f)_{L^2(\lambda)} - (\phi' \circ u, (-\Delta)^{-\beta_1} D_2 f)_{L^2(\lambda)}$$
$$- (v, (-\Delta)^{-\beta_1+\alpha_1} D_1 f)_{L^2(\lambda)}.$$

Now we are able to apply Corollary 7.1. I.e., assuming

$$\beta_2 \leq 2\beta_1,$$

suffices to obtain essential m-dissipativity in $L^2(W, \mu^\Phi, \mathbb{R})$ of the operator $(L_\Phi, \mathcal{F}C_b^\infty(B_W))$. To construct the corresponding diffusion process $X$ with infinite lifetime and providing a weak solution to the degenerate stochastic reaction–diffusion equation we need

$$\beta_2 > \frac{1}{2},$$

by (11) and

$$\alpha_1 = \alpha_2 \text{ or } \beta_1 + \frac{\alpha_2}{2} - \frac{\alpha_1}{2}, \beta_1 + \frac{\alpha_1}{2} - \frac{\alpha_2}{2} > \frac{1}{2},$$

by (13). In order to show hypocoercivity of the semigroup generated by $(L_\Phi, D(L_\Phi))$ and $L^2$-exponential ergodicity of $X$ we need to assume that $\sup_{t \in \mathbb{R}} |\phi'(t)| < \frac{\pi}{2}$ and

$$\beta_2 - \alpha_2 \leq 0 \text{ and } 2\beta_1 - \alpha_1 \leq 2\beta_2 - \alpha_2,$$

by (15).

**Acknowledgement** The first author gratefully acknowledges financial support in the form of a fellowship from the "Studienstiftung des deutschen Volkes".

**Funding** Open Access funding enabled and organized by Projekt DEAL. Not applicable.

**Availability of data and materials** Not applicable.

## Declarations

**Conflict of interest** The authors have no relevant financial or non-financial interests to disclose.

**Ethics approval** Not applicable.

**Code availability** Not applicable.

**Consent to participate** Not applicable.

**Consent for publication** Not applicable.

## References

1. Angiuli, L., Ferrari, S., Pallara, D.: Gradient estimates for perturbed Ornstein–Uhlenbeck semigroups on infinite-dimensional convex domains. J. Evol. Equ. **19**, 677–715 (2019). https://doi.org/10.1007/s00028-019-00491-y
2. Baudoin, F., Gordina, M., Herzog, D.P.: Gamma calculus beyond Villani and explicit convergence estimates for Langevin dynamics with singular potentials. Arch. Ration. Mech. Anal. **241**, 765–804 (2021). https://doi.org/10.1007/s00205-021-01664-1
3. Bertram, A., Grothaus, M.: Essential m-dissipativity and hypocoercivity of Langevin dynamics with multiplicative noise. J. Evol. Equ. **22**, 11 (2022). https://doi.org/10.1007/s00028-022-00773-y
4. Beznea, L., Boboc, N., Röckner, M.: Markov processes associated with $L^p$-resolvents and applications to stochastic differential equations on Hilbert spaces. J. Evol. Equ. **6**(4), 745–772 (2006) https://doi.org/10.1007/s00028-006-0287-2
5. Camrud, E., Herzog, D.P., Stoltz, G., Gordina, M.: Weighted $L^2$-contractivity of Langevin dynamics with singular potentials. Nonlinearity **35**(2), 998 (2021). https://doi.org/10.1088/1361-6544/ac4152
6. Conrad, F., Grothaus, M.: Construction, ergodicity and rate of convergence of N-particle Langevin dynamics with singular potentials. J. Evol. Equ. **10**(3), 623–662 (2010). https://doi.org/10.1007/s00028-010-0064-0
7. Conrad, F.: Construction and analysis of Langevin dynamics in continuous particle systems. Ph.D. thesis, Technische Universität Kaiserslautern, Published by Verlag Dr. Hut, München (2011)
8. Da Prato, G.: An Introduction to Infinite-Dimensional Analysis. Springer-Verlag, Berlin, Heidelberg (2006)
9. Da Prato, G., Lunardi, A.: Sobolev regularity for a class of second order elliptic PDE's in infinite dimension. Ann. Probab. **42**(5), 2113–2160 (2014). https://doi.org/10.1214/14-AOP936
10. Da Prato, G., Tubaro, L.: Self-adjointness of some infinite-dimensional elliptic operators and application to stochastic quantization. Probab. Theory Relat. Fields **118**, 131–145 (2000). https://doi.org/10.1007/s004400000073
11. Da Prato, G., Zabczyk, J.: Second order partial differential equations in Hilbert spaces. In: Mathematical Society Lecture Notes, vol. 293. Cambridge University Press, London (2002)
12. Da Prato, G.: Applications croissantes et equations d'evolution dans les espaces Banach. Academic Press, London (1976)
13. Dolbeault, J., Mouhot, C., Schmeiser, C.: Hypocoercivity for linear kinetic equations conserving mass. Trans. Am. Math. Soc. **367**(6), 3807–3828 (2015). https://doi.org/10.1090/S0002-9947-2015-06012-7
14. Eberle, A., Guillin, A., Zimmer, R.: Couplings and quantitative contraction rates for Langevin dynamics. Ann. Probab. **47**(4), 1982–2010 (2019). https://doi.org/10.1214/18-AOP1299
15. Eberle, A.: Uniqueness and Non-uniqueness of Semigroups Generated by Singular Diffusion Operators. Springer, Berlin (1999)
16. Eisenhuth, B., Grothaus, M.: Essential m-dissipativity for possibly degenerate generators of infinite-dimensional diffusion processes. Integr. Equ. Oper. Theory **94**, 28 (2022). https://doi.org/10.1007/s00020-022-02707-2
17. Grothaus, M., Stilgenbauer, P.: Hypocoercivity for Kolmogorov backward evolution equations and applications. J. Funct. Anal. **267**(10), 3515–3556 (2014). https://doi.org/10.1016/j.jfa.2014.08.019
18. Grothaus, M., Stilgenbauer, P.: Hilbert space hypocoercivity for the Langevin dynamics revisited. Methods Funct. Anal. Topol. **22**(2), 152–168 (2016)

19. Grothaus, M., Stilgenbauer, P.: A hypocoercivity related ergodicity method for singularly distorted non-symmetric diffusions. Integr. Equ. Oper. Theory **83**, 331–379 (2015). https://doi.org/10.1007/s00020-015-2254-1
20. Grothaus, M., Wang, F.-Y.: Weak poincaré inequalities for convergence rate of degenerate diffusion processes. Ann. Probab. **47**(5), 2930–2952 (2019). https://doi.org/10.1214/18-AOP1328
21. Herzog, D.P., Mattingly, J.C.: Ergodicity and Lyapunov functions for Langevin dynamics with singular potentials. Commun. Pure Appl. Math. **72**(10), 2231–2255 (2019). https://doi.org/10.1002/cpa.21862
22. Isserlis, L.: On a formula for the product-moment coefficient of any order of a normal frequency distribution in any number of variables. Biometrika **12**, 134–139 (1918). https://doi.org/10.2307/2331932
23. Ma, Z.-M., Röckner, M.: Introduction to the Theory of (Non-symmetric) Dirichlet Forms. Springer, Berlin (1992)
24. Mertin, M., Grothaus, M.: Hypocoercivity of Langevin-type dynamics on abstract smooth manifolds. Stoch. Process. Their Appl. **146**, 22–59 (2022). https://doi.org/10.1016/j.spa.2021.12.007
25. Pedersen, G.K.: Analysis Now. Graduate Texts in Mathematics, vol. 118. Springer, New York (1989)
26. Prévôt, C., Röckner, M.: A Concise Course on Stochastic Partial Differential Equations. Springer, Berlin (2007)
27. Villani, C.: Hypocoercivity. Mem. Am. Math. Soc. **202**(950), 141 (2009). https://doi.org/10.1090/S0065-9266-09-00567-5
28. Wang, F.-Y.: Hypercontractivity and applications for stochastic Hamiltonian systems. J. Funct. Anal. **272**(12), 5360–5383 (2017). https://doi.org/10.1016/j.jfa.2017.03.015